\documentclass[11pt,reqno]{elsarticle}

\usepackage{placeins}

\makeatletter
\def\ps@pprintTitle{%
 \let\@oddhead\@empty
 \let\@evenhead\@empty
 \def\@oddfoot{}%
 \let\@evenfoot\@oddfoot}
\makeatother

\usepackage{amsmath,amsthm,amscd,wrapfig,amsfonts,amssymb,mathtools}
\usepackage{enumitem,comment,mathabx,mathrsfs}
\usepackage{booktabs,multirow,lscape,cool,bm}
\usepackage{url,parskip,caption}
\usepackage[top=1.0in, bottom=1.0in, left=1.0in, right=1.0in]{geometry}

\newtheorem{theorem}{Theorem}

\newtheorem{remark}{Remark}

\newtheorem{conjecture}{Conjecture}

\providecommand{\e}[1]{\ensuremath{\times 10^{#1}}}

\makeatletter
\g@addto@macro\normalsize{%
  \setlength\abovedisplayskip{.4em}
  \setlength\belowdisplayskip{.4em}
  \setlength\abovedisplayshortskip{.4em}
  \setlength\belowdisplayshortskip{.4em}
}

\providecommand{\PS}[1]{\ensuremath{\mathit{Ps}\hskip.05em{}_{#1}}}

\providecommand{\PsiP}[1]{\ensuremath{\mathit{\Psi S}\hskip.05em{}_{#1}}}

\begin{document}

\begin{frontmatter}


\title
{
On the numerical evaluation of the prolate spheroidal wave functions of order zero
}

\begin{abstract}
We describe a method for the numerical evaluation of the angular prolate spheroidal wave functions of the first kind
of order zero.
It is based on the observation that underlies the WKB method, namely that many second
order differential equations admit solutions whose logarithms can be represented much
more efficiently than the solutions themselves.  However, rather than exploiting this fact
to construct asymptotic expansions of the prolate spheroidal wave functions, 
our algorithm operates by  numerically solving the Riccati equation satisfied 
by their logarithms.      Its running time grows  much more slowly with bandlimit
and characteristic exponent than standard algorithms.  
We illustrate  this and other properties of our algorithm  with numerical experiments.

\end{abstract}

\author{James Bremer}
\ead{bremer@math.toronto.edu}
\address{Department of Mathematics, University of Toronto}

\end{frontmatter}
\begin{section}{Introduction}

Many families of special functions are defined by second order linear ordinary differential equations
whose coefficients depend on one or more parameters.  The cost to represent solutions
of such equations using standard methods, such as polynomial and trigonometric expansions,
typically grows quite quickly with magnitudes of the parameters.  
It is well known, though, that in many cases there exist solutions whose logarithms can be represented 
at a cost which is bounded independent of the parameters, or at least grows extremely slowly with them.

This observation is the basis of  many approaches to the asymptotic approximation of the solutions of  
second order linear ordinary differential equations.  
The WKB method, for instance, can be used to construct asymptotic expansions
of the solutions of equations of the form
\begin{equation}
y''(t) + \lambda^2 q(t) y(t) = 0
\label{introduction:ode}
\end{equation}
with $q$ a strictly positive smooth function.  In this case, 
there exist a solution $y$ of (\ref{introduction:ode}) and 
a sequence of functions $r_0,r_1,\ldots$  
which depend on $q$ and its derivatives, but not $\lambda$, such that 
\begin{equation}
y(t) = \exp\left( i\sum_{n=0}^N \lambda^{1-n}\ r_n(t)\right)
\left(1 + \mathcal{O}\left(\lambda^{-N}\right)\right)
\ \ \mbox{as} \ \ \lambda \to \infty
\label{introduction:asym}
\end{equation}
(see, for instance, Section~7.2 of \cite{Miller}).
The solution $y$ is oscillatory and the  cost of representing it with standard methods grows 
linearly with $\lambda$.    By contrast, the functions $r_n$ can be represented at a 
cost which depends only on the complexity of $q$, and not on $\lambda$.  This means that for sufficiently large values of $\lambda$,
it is much more efficient to represent $y$ via (\ref{introduction:asym}) than to construct a polynomial expansion of $y$ itself.


Asymptotic methods allow for the evaluation of many special functions
in time independent of the parameters they depend on,
but are only  accurate for sufficiently large values of those parameters.
They are often coupled with methods which are accurate at all values of the  parameters, 
but whose running times grow with the parameters, to obtain efficient
algorithms for the numerical evaluation of a family of special functions.
This is the approach used in   \cite{Bogaert-Michiels-Fostier}, for instance,
to evaluate the Legendre  polynomials in time independent of their degree.  
Further examples can be found in \cite{DLMF,Olver,Gil-Segura-Temme} and their references.
Unfortunately, in many cases existing asymptotic expansions are either not amenable
to numerical evaluation or only achieve high-accuracy at extremely large values of the parameters,
making such an approach infeasible.  

The prolate spheroidal wave functions of order zero and bandlimit $\gamma > 0$
are an example of a family of special functions for which such an approach is not
viable.  They are the solutions of the second order linear ordinary differential equation
\begin{equation}
(1-z^2)y''(z) -2z y'(z) + (\chi-\gamma^2z^2) y(z) = 0,
\label{introduction:rswe}
\end{equation}
which we call the reduced spheroidal wave equation (reduced because it is obtained from the 
more familiar spheroidal wave equation by deleting one parameter).
The solutions of most interest are those which satisfy the boundary conditions
\begin{equation}
\lim_{z \to  \pm 1} y'(z) \sqrt{1-z^2}=0.
\label{introduction:bc}
\end{equation}
Together  (\ref{introduction:rswe}) and (\ref{introduction:bc}) constitute a
singular self-adjoint Sturm-Liouville problem.  Consequently, there exists a sequence
\begin{equation}
\chi_0(\gamma) < \chi_1(\gamma) < \chi_2(\gamma) < \cdots
\label{introduction:chi}
\end{equation}
of values of the parameter $\chi$ for which solutions of (\ref{introduction:rswe}) 
satisfying (\ref{introduction:bc}) exist.   For each $\chi_n(\gamma)$, there is a corresponding one-dimensional space
of solutions of (\ref{introduction:rswe}), and we use $\PS{n}(z;\gamma)$ to denote
the particular element of that subspace which agrees with the Legendre function $P_n(z)$ at the point $z=0$.
The parameter $n$ usually called  the characteristic exponent ---  this  term comes from the standard
mechanism used to define $\PS{\nu}(z;\gamma)$ for noninteger $\nu$   ---
and $\PS{n}(z;\gamma)$ is known as the angular prolate spheroidal wave function of the first kind
of bandlimit $\gamma$, order $0$ and  characteristic exponent $n$.

Although uniform asymptotic expansions of the $\PS{n}(z;\gamma)$ are available \cite{Dunster},  they 
only achieve high accuracy for extremely large values of $\gamma$  and 
involve a complicated change of variables,  thus making them difficult to exploit in numerical computations.
This leaves algorithms whose running times increase fairly rapidly with $\gamma$ and $n$ as the only viable
mechanisms for the numerical evaluation of $\PS{n}(z;\gamma)$.
The standard approach  is the  Osipov-Xiao-Rokhlin algorithm   \cite{Osipov-Rokhlin,Osipov-Rokhlin-Xiao}.
It proceeds by solving an eigenproblem for a symmetric tridiagonal matrix in order to construct a Legendre expansion
representing the desired prolate function. The dependence of the running time  of this algorithm
on the parameters $\gamma$ and $n$ is not fully understood, but the numerical experiments of \cite{Feichtinger}
suggest that it behaves as  $\mathcal{O}\left( n + \sqrt{\gamma n}\right)$
for large values of the parameters.

We describe an algorithm for evaluating $\PS{n}(z;\gamma)$
whose running time grows much more slowly with $\gamma$ and $n$.
It operates by solving the Riccati equation satisfied by the logarithms
of the solutions of the reduced spheroidal wave equation numerically. 
Most solutions of that equation are no easier to represent than the prolate spheroidal 
wave function themselves.  However, standard results (such as those appearing in Section~7.2 of \cite{Miller})
imply that there is a  ``WKB solution'' of the reduced spheroidal wave equation whose logarithm
can be  asymptotically approximated by a series of relatively simple functions.
It is straightforward to identify this solution from well-known formulas, and it is 
the logarithm of this solution we construct by solving the Riccati equation.   Based on the numerical experiments discussed 
in this article, we believe that  the cost of representing it via polynomial expansions grows sublogarithmically 
with $\gamma$ and is bounded  independent of $n$ for fixed $\gamma$.    The running time of our algorithm displays the same dependence
on $\gamma$ and $n$.

We also state several conjectures about the monotonicity properties of the modulus
of this WKB solution.  These conjectures are relevant  because the logarithm of any solution
of a second order differential equation is related to the modulus of the solution
through a simple formula.  Our conjectures assert that the monotonicity properties of the reduced spheroidal
wave equation are similar to those of Legendre's differential equation, which it generalizes.

There are two significant limitations of the algorithm described in this paper.  First, it
requires knowledge of Sturm-Liouville eigenvalue $\chi_n(\gamma)$
in order to calculate $\PS{n}(z;\gamma)$.
Second, while the running time of our algorithm grows much more slower
with  $\gamma$ and $n$ than does the  Osipov-Rokhlin-Xiao algorithm, the later is more efficient until fairly large
values of $\gamma$ and $n$.  The author remedies these problems in a separate article, which describes
a method for computing the value of $\chi_n(\gamma)$ in time independent the parameters as well as 
an approach to  accelerating the algorithm of this paper which makes it faster than the Osipov-Xiao-Rokhlin
algorithm except at extremely small values of the parameters.   A comparison of the running times of the 
unaccelerated and accelerated versions of the algorithm of this paper and the Xiao-Osipov-Rokhlin method can be found 
in Section~\ref{section:experiments}  of this paper.

The remainder of this article is structured as follows.  Section~\ref{section:riccati} briefly
discuss the Riccati equation satisfied by the logarithms of solutions of second order linear ordinary
differential equations.  In Section~\ref{section:spheroidal}, we discuss the 
prolate spheroidal wave functions of order zero.
In Section~\ref{section:monotonicity}, the monotonicity properties of Legendre's differential
equation are described and  we make several conjectures regarding the monotonicity 
properties of the reduced spheroidal wave equation.
Section~\ref{section:algorithm} details our  numerical
algorithm for the evaluation of the prolate spheroidal wave functions of order zero.
 In Section~\ref{section:experiments}, we describe the results of numerical experiments
which were conducted to demonstrate the properties of our algorithm.  
We close with a few brief remarks in  Section~\ref{section:conclusion}.

\end{section}

\begin{section}{Riccati's equation and its variants}
\label{section:riccati}

A straightforward calculation shows that if  $y(x) = \exp(r(x))$ solves 
the second order differential equation
\begin{equation}
y''(x) + q(x) y(x) = 0,
\label{riccati:ode}
\end{equation}
then $r$ satisfies the Riccati equation
\begin{equation}
r''(x) + (r'(x))^2 + q(x) = 0.
\label{riccati:riccati}
\end{equation}
We note that under extremely mild regularity conditions on its coefficients, any second order linear ordinary differential equation 
can be put into the form (\ref{riccati:ode}) through a simple transformation (see, for instance, Section~5.6 of \cite{Hille}).

By inserting the expression $r(x) = \pm i \psi(x) + \beta(x)$  into (\ref{riccati:riccati}),
it can be shown that if  $\psi$ and $\beta$  satisfy the system of equations
\begin{equation}
\left\{
\begin{aligned}
\beta''(x) +(\beta'(x))^2-(\psi'(x))^2+q(x)&=0\\
 \psi''(x) + 2\psi'(x) \beta'(x) &= 0,
\end{aligned}
\right.
\label{riccati:system}
\end{equation}
then $r$ solves (\ref{riccati:riccati}).  The second equation in (\ref{riccati:system})
admits the formal solution 
\begin{equation}
\beta(x) = -\frac{1}{2}\log(\psi'(x)),
\label{riccati:beta}
\end{equation}
which when inserted  into the first equation in (\ref{riccati:system}) results in 
\begin{equation}
q(x) - (\psi'(x))^2 + \frac{3}{4} \left(\frac{\psi''(x)}{\psi'(x)}\right)^2
-\frac{1}{2}\frac{\psi'''(x)}{\psi'(x)}=0.
\label{riccati:kummer}
\end{equation}
Equation (\ref{riccati:kummer}) is known as  Kummer's equation after E.~E.~Kummer who studied it in \cite{Kummer}.
Clearly,  if $\psi$ does not vanish in an open interval and it satisfies Kummer's equation  there, then
\begin{equation}
r_1(x) = i \psi(x) - \frac{1}{2}\log(\psi'(x)) \ \ \ \mbox{and}\ \ 
r_2(x) = -i \psi(x) - \frac{1}{2}\log(\psi'(x))
\end{equation}
are solutions of Riccati's equation on that interval.  It follows that 
\begin{equation}
u(x) = \frac{\sin(\psi(x))}{\sqrt{\psi'(x)}} \ \ \mbox{and}\ \ 
v(x) = \frac{\cos(\psi(x))}{\sqrt{\psi'(x)}}
\label{riccati:uv}
\end{equation}
constitute a basis in the space of solutions of  (\ref{riccati:ode}).
In this event, we refer to $\psi$ as a phase function for (\ref{riccati:ode}).
We note that the realization of the square root used in (\ref{riccati:uv}) is immaterial.
For obvious reasons, 
\begin{equation}
m(x) = (u(x))^2 + (v(x))^2 = \frac{1}{\psi'(x)}
\label{riccati:modulus}
\end{equation}
is known as the modulus function associated with $\psi$.  It can be readily verified 
that $m$ satisfies the differential equation
\begin{equation}
m'''(x) + 4 q(x) m'(x) + 2 q'(x) m(x) = 0,
\label{riccati:appell}
\end{equation}
which we refer to  as Appell's equation after P.~E.~Appell, who discussed it in \cite{Appell}.

Given any pair $u, v$  of solutions of (\ref{riccati:ode}) whose Wronskian is $w \neq 0$
and such that  such that the modulus function (\ref{riccati:modulus}) is nonzero in an open
interval $I$, it can be shown by a straightforward calculation that
\begin{equation}
\psi'(x) = \frac{w}{(u(x))^2+(v(x))^2}
\label{riccati:alphap}
\end{equation}
satisfies (\ref{riccati:kummer}) on that interval.  It follows that any antiderivative of $\psi'$ 
is a phase function for (\ref{riccati:ode}) on $I$.  Requiring that (\ref{riccati:uv}) holds
determines $\psi$  up to an integral multiple of $2\pi$, but further restrictions
are required to determine it uniquely.
We note that on any set where $u$ and $v$ are real-valued, the modulus function
(\ref{riccati:modulus}) is nonzero since independent solutions of a second
order differential equation cannot simultaneously vanish.

Because of the close relationship between the functions $\psi$, $m$ and $r$, we regard
them all as representations of the logarithm of a solution of (\ref{riccati:ode}).
Moreover, we view Equations~(\ref{riccati:riccati}), (\ref{riccati:kummer}) and 
(\ref{riccati:appell}) as essentially interchangeable mechanisms for
computing the logarithm of a solution of (\ref{riccati:ode}).

\end{section}

\begin{section}{The prolate spheroidal wave functions of order zero}
\label{section:spheroidal}

In this section, we briefly discuss certain facts regarding the 
prolate spheroidal wave functions of order zero which will be 
used in the algorithm of this paper.

\begin{subsection}{The angular prolate spheroidal wave functions of the first kind}

In addition to being the solutions of a singular self-adjoint Sturm-Liouville,
the angular prolate spheroidal wave functions of bandlimit $\gamma > 0$, order zero and integer characteristic exponents
\begin{equation}
\PS{0}(z;\gamma),\ \PS{1}(z;\gamma),\ \PS{2}(z;\gamma),\ \ldots
\label{spheroidal:ps}
\end{equation}
%
are the eigenfunctions of the restricted Fourier operator
\begin{equation}
\mathscr{F}_\gamma\left[f\right](z) = \int_{-1}^1 \exp\left(i\gamma z t \right) f(t)\ dt.
\label{spheroidal:fourier}
\end{equation}
As such, they provide an efficient mechanism for representing elements of the image of $\mathscr{F}_\gamma$,
which is the set of functions with bandlimit $\gamma$.
The dual nature of the functions (\ref{spheroidal:ps}) was widely publicized in an article \cite{PSWFI} published in the 1960s,
but it was known much earlier (see, for instance, Section~3.8 of \cite{Meixner} and the references cited there).

It is shown in \cite{Landau-Widom}  that the magnitudes of the first $2/\pi \gamma$ or so eigenvalues of $\mathscr{F}_\gamma$
are close to  $\sqrt{2\pi/\gamma}$, the magnitudes of the next $\mathcal{O}\left(\log(\gamma)\right)$ eigenvalues
decay extremely rapidly, and the remaining eigenvalues are all close
to zero.  It follows that
 only the first 
$2/\pi \gamma + \mathcal{O}\left(\log(\gamma)\right)$ of the functions (\ref{spheroidal:ps})
are needed to represent elements of the image of $\mathscr{F}_\gamma$, which is the
space of functions with bandlimit $\gamma$,  with extremely high relative accuracy.
In other words, for the purposes of numerical computation, the dimension
of the space of functions with bandlimit $\gamma$ is  $2/\pi \gamma + \mathcal{O}\left(\log(\gamma)\right)$.

The reduced spheroidal wave equation (\ref{introduction:rswe}) has a regular singular point at $z=1$
and zero is a double root of the corresponding indicial equation.
Consequently, there is a one-dimensional subspace of solutions which are regular at $z=1$ and all  other
solutions have logarithmic singularities there (see, for instance, Chapter~5 of \cite{Hille}).  
Because it is an eigenfunction of the restricted Fourier operator,  $\PS{n}(z;\gamma)$ 
is entire and so the subspace of solutions which are regular at $z=1$ comprises its multiples.  Every other solution has a 
logarithmic singularity at $1$.

Since $\PS{n}(z;\gamma)$ is an eigenfunction
of (\ref{spheroidal:fourier}), there exists a constant $A_n(\gamma)$ such that 
\begin{equation}
\PS{n}(z;\gamma)  =
A_n(\gamma) \frac{\sin\left(\gamma z\right) }{\gamma z} 
\left(1 + \mathcal{O}\left(\frac{1}{z}\right) \right)\ \ \mbox{as} \ \ z \to \infty.
\label{spheroidal:psasym}
\end{equation}
It is a consequence of this and the fact that $\PS{n}(z;\gamma)$ is regular at $1$ that 
$\PS{n}(z;\gamma)$ is an element of $L^2(1,\infty)$.

\end{subsection}

\begin{subsection}{The radial prolate spheroidal wave function of the third kind of order zero}
Another solution of the reduced spheroidal wave equation, which is known as the radial prolate spheroidal wave function
of the third kind of order zero, is defined via the formula
\begin{equation}
S^{(3)}_n(z;\gamma) =  
\frac{1}{\PS{n}(1,\gamma)}
\int_1^\infty \exp(i\gamma zt) \PS{n}(t; \gamma)\ dt.
\label{spheroidal:s3}
\end{equation}
Since $\PS{n}(z;\gamma)$ is in $L^2(1,\infty)$, 
the integral is  absolutely convergent for all $\mbox{Im}(z) > 0$ and
 $S^{(3)}_n(z; \gamma)$ is usually taken to be function obtained by analytically continuing
it to the cut plane $\mathbb{C}\setminus\{-\infty,1\}$.  However, for our purposes, it is more convenient
to take its domain to be  $\{ z : \mbox{Im}(z) \geq 0 \} \setminus \{-1,1\}$.
The asymptotic behaviour of $S^{(3)}(z;\gamma)$ is obvious from Formula~(\ref{spheroidal:s3}):
\begin{equation}
S^{(3)}_n(z;\gamma)  =
\frac{\exp\left(i\gamma z \right)}
{\gamma z}
\left(1+ \mathcal{O}\left(\frac{1}{z}\right)\right)\ \ \ \mbox{as} \ \ z \to \infty.
\label{spheroidal:s3asym}
\end{equation}

\end{subsection}

\begin{subsection}{Two variants of the reduced spheroidal wave equation}

We shall make use of two variants of the reduced spheroidal wave equation of the form (\ref{riccati:ode}).
The first of these is 
\begin{equation}
y_1''(x) + \left( 
\frac{\exp(x)-\frac{1}{4}}{\left(1-2\exp(x)\right)^2}
+ \chi\ \frac{2\exp(x) - 1}{(1-2\exp(x))^2} 
-\gamma^2\exp(-2x) 
\right)y_1(x) = 0, \ \ \ 0 < x < \infty,
\label{spheroidal:erswe}
\end{equation}
which is satisfied by any function of the form $y_1(t) = y\left(1-\exp(-x)\right)\sqrt{2-\exp(-x)}$ 
with $y$ a solution of (\ref{introduction:rswe}).  We refer to (\ref{spheroidal:erswe}) as the 
exponential form of the reduced spheroidal wave equation.  Our motivation for introducing
this exponential change of variables will be explained in Section~\ref{section:algorithm}.

Similarly, if $y$ is a solution of (\ref{introduction:rswe}) on the upper half
of the imaginary axis, then  $y_2(t) = y(it)\sqrt{1+t^2}$ satisfies 
\begin{equation}
y_2''(t) -  \left(\frac{1}{\left(1+t^2\right)^2}
+ \frac{\chi+\gamma^2 t^2}{1+t^2}
\right)y_2(t) = 0,\ \ \ 0 < t < \infty.
\label{spheroidal:irswe}
\end{equation}

\end{subsection}

\begin{subsection}{The phase and modulus functions associated with $S^{(3)}_n(z;\gamma)$}

Since the coefficient in (\ref{spheroidal:erswe}) is real-valued,
the real and imaginary parts of 
\begin{equation}
S^{(3)}_n(1-\exp(-x);\gamma)\sqrt{2-\exp(-x)}
\end{equation}
are separately solutions. 
Accordingly, 
\begin{equation}
M_n(1-\exp(-x);\gamma) (2-\exp(-x)),
\end{equation}
where 
\begin{equation}
M_n(z;\gamma) = 
\left|S^{(3)}_n(z;\gamma)\right|^2,
\end{equation}
is a modulus function for (\ref{spheroidal:erswe}).     
Since the Wronskian of any pair of solutions
of (\ref{spheroidal:erswe}) is necessarily constant, it follows from (\ref{spheroidal:s3asym})
that the Wronskian of the pair of solutions consisting of the 
real and imaginary parts of $S^{(3)}_n(z;\gamma)$ is $\gamma$
and we define $\PsiP{n}(x;\gamma)$ on the interval $[0,\infty)$ via 
\begin{equation}
\PsiP{n}(x;\gamma)  = \int_{\infty}^x \frac{\gamma}{M_n(1-\exp(-u);\gamma)(2-\exp(-u))}\ du.
\label{spheroidal:psiprime}
\end{equation}
It follows from the  discussion in Section~\ref{section:riccati} that $\PsiP{n}(x;\gamma)$ is a phase
function for the exponential form of the reduced spheroidal wave equation.  In particular,
\begin{equation}
\frac{\sin\left(\PsiP{n}(x;\gamma)\right)}{\sqrt{\frac{d\PsiP{n}}{dx} (x;\gamma) }}\ \ \mbox{and} \ \ \
\frac{\cos\left(\PsiP{n}(x;\gamma)\right)}{\sqrt{\frac{d\PsiP{n}}{dx} (x;\gamma) }}
\label{spheroidal:phasebasis}
\end{equation}
form a basis in the space of solutions of Equation~(\ref{spheroidal:erswe}).    
Moreover, (\ref{spheroidal:psiprime}) ensures that 
\begin{equation}
\lim_{x\to \infty} \PsiP{n}(x;\gamma)=0.
\label{spheroidal:psipconstant}
\end{equation}
Since $\PS{n}(1-\exp(-x);\gamma)\sqrt{2-\exp(-x)},$ is a solution of the exponential form
of the reduced spheroidal wave equation, there exist constants $C_n(\gamma)$ and $D_n(\gamma)$ such that
\begin{equation}
\PS{n}(1-\exp(-x);\gamma)\sqrt{2-\exp(-x)}
=
C_n(\gamma) \frac{\sin\left(\PsiP{n}(x;\gamma)\right)}{\sqrt{\frac{d\PsiP{n}}{dx} (x;\gamma) }}
+ D_n(\gamma) \frac{\cos\left(\PsiP{n}(x;\gamma)\right)}{\sqrt{\frac{d\PsiP{n}}{dx} (x;\gamma) }}.
\label{spheroidal:pscombo1}
\end{equation}
 In fact, we claim that
condition (\ref{spheroidal:psipconstant}) ensures $D_n(\gamma)=0$ in (\ref{spheroidal:pscombo1}); that is,
our choice of the constant of integration for $\PsiP{n}(x;\gamma)$ guarantees that 
$\PS{n}(1-\exp(-x);\gamma)\sqrt{2-\exp(-x)}$  is a multiple of 
$$\frac{\sin\left(\PsiP{n}(x;\gamma)\right)}{\sqrt{\frac{d\PsiP{n}}{dx} (x;\gamma) }}.$$
To see this, we first observe that (\ref{spheroidal:pscombo1})  holds if and only if 
if
\begin{equation}
\begin{aligned}
\PS{n}(1-\exp(-x);\gamma)
&=
\frac{C_n(\gamma)}{\sqrt{\gamma}} \sqrt{M_n(1-\exp(-x);\gamma)}\sin\left(\PsiP{n}(x;\gamma)\right)\\
&+ \frac{D_n(\gamma)}{\sqrt{\gamma}} \sqrt{M_n(1-\exp(-x);\gamma)}\cos\left(\PsiP{n}(x;\gamma)\right).
\end{aligned}
\label{spheroidal:pscombo2}
\end{equation}
Next, we note that  since  every solution of the reduced spheroidal wave equation which is
not a multiple of $\PS{n}(z;\gamma)$ has a logarithmic singularity at $1$,
\begin{equation}
\lim_{x\to \infty} M_n(1-\exp(-x);\gamma) = \lim_{z\to 1^-} M_n(z;\gamma) = \infty.
\label{spheroidal:limit}
\end{equation}
Equations~ (\ref{spheroidal:psipconstant})  and (\ref{spheroidal:limit}) imply that
\begin{equation}
\lim_{x\to \infty} 
\sqrt{M_n(1-\exp(-x);\gamma)}\cos\left(\PsiP{n}(x;\gamma)\right) = \infty,
\end{equation}
whereas $\PS{n}(x;\gamma)$ is nonsingular at $x=1$.  It follows that the constant $D_n(\gamma)$ 
in (\ref{spheroidal:pscombo2}) must be 0.


\end{subsection}

\end{section}

\begin{section}{Montonicity properties of the reduced spheroidal wave equation}
\label{section:monotonicity}

Many second order differential equations admit solutions whose logarithms
are easier to represent, either numerically or symbollically, then
the solutions themselves.  Because of the close relationship
between the logarithms of solutions of  second order linear ordinary differential equations
and  their modulus functions, this is often demonstrated
by establishing the existence of modulus functions satisfying certain  monotonicity properties.

From (\ref{spheroidal:s3asym}), it is clear that $S^{(3)}_n(z;\gamma)$ is a ``WKB solution''
of the reduced spheroidal wave equation in that, at least for large values of $z$, its logarithm 
can be approximated by a polynomial expansion at cost which is independent of $\gamma$ and $n$.  
In fact, plotting its modulus function for  various values of $n$ and $\gamma$ suggests that 
the modulus function associated with $S^{(3)}_n(z;\gamma)$ satisfies rather strong monotonicity properties
(see Figure~\ref{figure:plots} for representative plots).  

In this section, we first discuss the monotonicity properties of a certain  modulus function for  Legendre's differential
equation.  Then, we conjecture that the modulus function corresponding to $S_\nu(z;\gamma)$ behaves in a similar fashion.

\begin{figure}[!h]
\begin{center}
\hfil
\includegraphics[width=.45\textwidth]{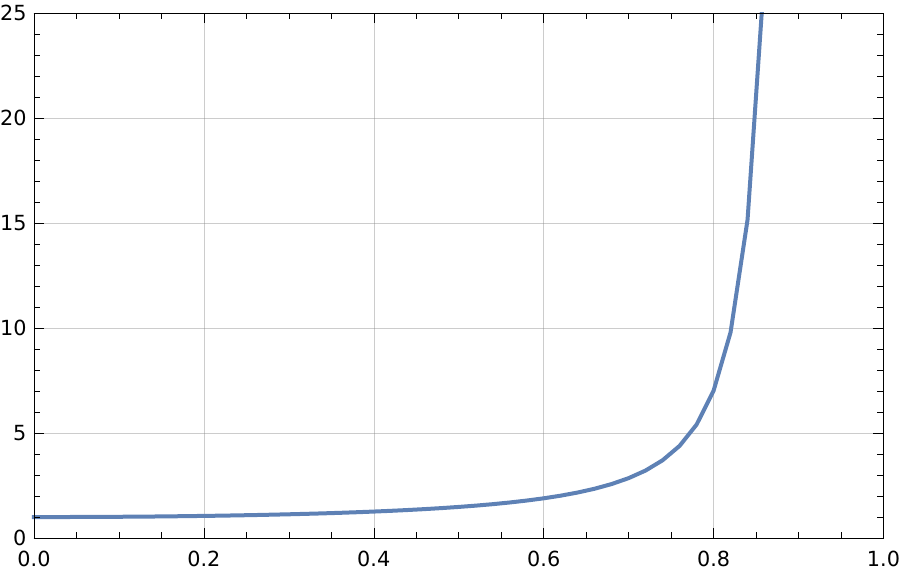}
\hfil
\includegraphics[width=.45\textwidth]{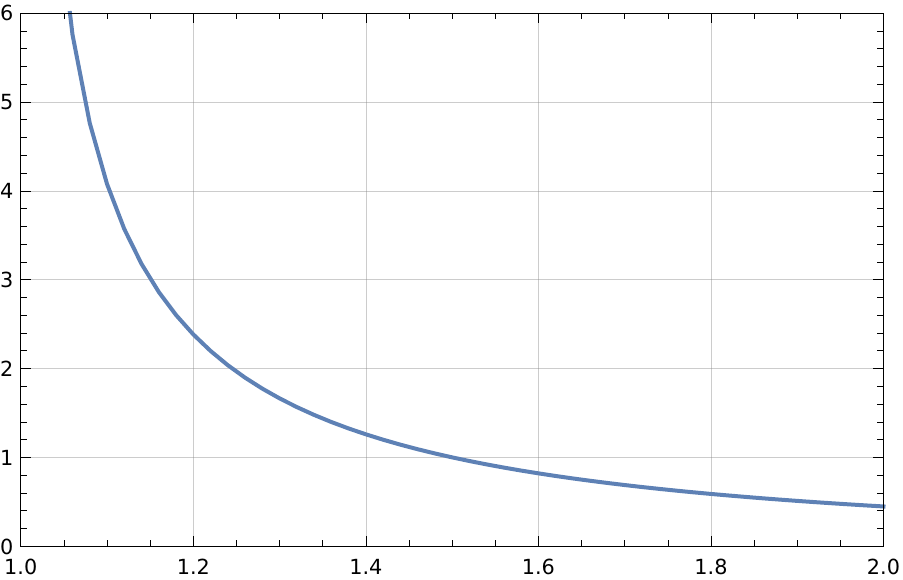}
\hfil
\end{center}
\caption{
On the left is a plot of 
 $\left|S^{(3)}_n(x;\gamma)\right|^2$ over the interval $(0,1)$ when $\gamma=20$ and $n=10$.  
On the right, is a plot of the same function over the interval $(1,2)$.
}
\label{figure:plots}
\end{figure}

\begin{subsection}{Completely, absolutely and multiply monotone functions}

A function $f$ defined on an open interval $I$ is $k$-times monotone there provided
$(-1)^j f^{(j)}(z) \geq 0$ for all nonnegative integers $j \leq k$ and all $z \in I$.
It is completely monotone if $(-1)^k f^{(k)}(z) \geq 0$
for all nonnegative integers $k$ and all $z\in I$.
Finally,  $f$ is absolutely monotone if $f^{(k)}(z) \geq 0$ for all nonnegative
$k$ and all $z\in I$.  
It is well known that $f$ is completely monotone on $(0,\infty)$ if and only if 
it is the Laplace transform of a nonnegative Borel measure
(see, for instance,  Chapter~4 of \cite{Widder} for a proof of this).
\end{subsection}

\begin{subsection}{Legendre's differential equation}

Legendre's differential equation, which is the  special case of  (\ref{introduction:rswe}) when  $\gamma=0$, 
is a classical example of a differential equation admitting a modulus function which satisfies
various monotonicity properties.
The Legendre function of the second kind of degree $\nu$ is given by
\begin{equation}
Q_\nu(z) =  \exp\left(-i \frac{\pi}{2}(\nu+1)\right)\int_0^\infty \exp(izt) j_\nu(t)\ dt,
\label{monotonicity:qnu}
\end{equation}
where $j_\nu$ denotes the spherical Bessel function of order $\nu$.
The integral is absolutely convergent for $\mbox{Im}(z) > 0$
and  $Q_\nu(z)$ is typically defined to be the function obtained by analytically
continuing the integral (\ref{monotonicity:qnu}) to the cut plane $\mathbb{C}\setminus\left(-\infty,1\right]$.
However, we prefer to takes its domain to be the set $\{z  : \mbox{Im}(z) \geq 0\}\setminus\{-1,1\}$.

It can be easily verified that $Q_\nu(x)\sqrt{1-x^2}$ satisfies the differential equation
\begin{equation}
y''(x) + \left(
\frac{1}{(1-x^2)^2} + \frac{\chi}{1-x^2}
\right)y(x) = 0,\ \ \ -1 < x < 1, 
\label{monotonicity:lege}
\end{equation}
where $\chi = \nu (\nu+1)$.
Since the coefficient in (\ref{monotonicity:lege})
is real-valued, the real and imaginary parts of this function are separately solutions,
and $\left|Q_\nu(x)\right|^2(1-x^2)$ is a modulus function for (\ref{monotonicity:lege}).
Because
\begin{equation}
j_\nu(t) \sim \frac{\sin\left(t-\frac{\pi}{2} \nu\right)}{t} \ \ \mbox{as} \ \ t \to \infty,
\end{equation}
the contour in (\ref{monotonicity:qnu}) can be shifted from the real axis to the imaginary axis
whenever $\mbox{Re}(z) > 1$.    This results in the formula
\begin{equation}
Q_\nu(z) =  \exp\left(-i \frac{\pi}{2}(\nu+1)\right) \int_0^\infty \exp(-zt) i_\nu(t)\ dt
\label{monotonicity:qnu2}
\end{equation}
with $i_\nu$ the modified spherical Bessel function of the first kind of order $\nu$.
Since $i_\nu$ is positive on $(0,\infty)$,  (\ref{monotonicity:qnu2}) implies that 
$Q_\nu(x)$ is completely monotone on the interval $(1,\infty)$.
From this and the formula
\begin{equation}
\left|Q_\nu(x)\right|^2 = 
 \int_1^\infty
Q_\nu\left(x^2+(1-x^2)t \right) \frac{dt}{\sqrt{t^2-1}},
\label{monotonicity:legendre}
\end{equation}
which holds for $\nu \geq 0$ and can be found in \cite{durand78}, it follows that  $\left|Q_\nu(x)\right|^2$ is absolutely monotone 
on the interval $(-1,1)$ when $\nu \geq 0$.  Since the square of a completely monotone function is completely monotone,
we have the following theorem summarizing the monotonicity properties of 
$\left|Q_\nu(x)\right|^2$:

\begin{theorem}
For fixed $\nu \geq 0$, $\left|Q_\nu(x)\right|^2$
is absolutely monotone on $(0,1)$ and completely monotone on $(1,\infty)$.
\end{theorem}

\end{subsection}

\begin{subsection}{The reduced spheroidal wave equation}

Many other second order linear ordinary differential equations admit modulus functions 
which satisfy various monotonicity properties.   Relevant formulas 
for the Jacobi functions, Gegenbauer functions and Hermite functions  
can be found in \cite{durand78},
and the articles \cite{hartman61} and \cite{hartman73} 
give conditions under which  a second order linear ordinary differential equation
admits a completely monotone modulus function.

However, to the author's knowledge, no results concerning the monotonicity properties
of the reduced spheroidal wave equation appear in the literature.
Since the reduced spheroidal wave equation (\ref{introduction:rswe}) generalizes Legendre's differential
equation, it is reasonable to suspect that the solution $S^{(3)}_n(z;\gamma)$
behaves similarly to $Q_n(z)$.  This is  further suggested by the 
many similaries between Formulas~(\ref{spheroidal:s3}) and (\ref{monotonicity:qnu}).
Accordingly, we make the following conjecture regarding the radial prolate
spheroidal wave function of the third kind:

\begin{conjecture}
For fixed $\gamma > 0$ and $n \geq 0$,  $M_n(z;\gamma)$
is absolutely monotone on $(0,1)$ and completely monotone on $(1,\infty)$.
\label{conjecture1}
\end{conjecture}

Since the coefficient in  (\ref{spheroidal:irswe}) is negative for $u \in (0,\infty)$,  the function $S_n^{(3)}(iu; \gamma)$ is nonoscillatory
for $0 < u < \infty$.    Based on experiments performed using computer algebra systems, we make the following conjecture
regarding the behavior of $S^{(3)}_n(z;\gamma)$ on the imaginary axis:

\begin{conjecture}
For fixed $\gamma >0$ and $n \geq 0$,  
$S_n^{(3)}(iu;\gamma)$  is $\left(2+n\right)$-times monontone on the interval $(0,\infty)$.
\label{conjecture2}
\end{conjecture}

\vskip 1em
\begin{remark}
In \cite{prolates2}, we discuss the standard mechanism for defining 
$S^{(3)}_\nu(z;\gamma)$ for noninteger values of $\nu$ and generalize
Conjectures~\ref{conjecture1} and \ref{conjecture2} to that case.
\end{remark}

\end{subsection}

\end{section}

\begin{section}{Numerical algorithm}
\label{section:algorithm}

We now describe our algorithm for the numerical evaluation of $\PS{n}(z;\gamma)$.
In addition to the values of $\gamma$ and $n$, it takes as input the Sturm-Liouville eigenvalue 
$\chi_n(\gamma)$ and a positive integer $k$.  The algorithm operates by 
constructing a $k^{th}$ order piecewise Chebyshev expansion
representing the phase function $\PsiP{n}(x;\gamma)$ on an interval of the form $\left[0,\beta\right)$.
More explicitly, the piecewise Chebyshev expansion comprises a partition
\begin{equation}
0 = x_0 < x_1 < x_2 < \ldots < x_n = \beta
\label{algorithm:part}
\end{equation}
of $\left[0,\beta\right)$ together with the coefficients $\{ a_{ij} : 1 \leq i \leq n,  \ 1 \leq j \leq k\}$
in the expansion
\begin{equation}
\PsiP{n}(x;\gamma) \approx \sum_{i=1}^{n} I_{\left[x_{i-1},x_{i}\right)}(x) \sum_{j=0}^k a_{ij}\ T_j\left(\frac{2}{x_{i}-x_{i-1}} x + 
\frac{x_{i-1}+x_{i}}{x_{i-1}-x_{i}}\right),
\label{algorithm:chebexp}
\end{equation}
where $T_j$ denotes the Chebyshev polynomial of degree and $I_{\left[x_{i-1},x_{i}\right)}(x)$ is the characteristic function 
\begin{equation}
I_{\left[x_{i-1},x_{i}\right)}(x) = 
\begin{cases}
1 & \mbox{if} \ \ x_{i-1} \leq x < x_i\\
0 & \mbox{otherwise}.
\end{cases}
\end{equation}
In all of the experiments described in this paper, we took $k=29$ and $\beta=10^{120}$.   
We will shortly explain why such a large value of $\beta$ is needed.

The phase function is related to $\PS{n}(z;\gamma)$ via the formula
\begin{equation}
\PS{n}(z;\gamma) = C_n(\gamma)  \frac{\sin\left(\PsiP{n}(x;\gamma)\right)}{\sqrt{(1+z)\frac{d\PsiP{n}}{dx}(x;\gamma) }},
\label{algorithm:psn}
\end{equation}
where $x=-\log(1-z)$ and $C_n(\gamma)$ is a constant which must be calculated.  Since we normalize $\PS{n}(z;\gamma)$ by requiring
that it equal the Legendre function $P_n(z)$ at the point $z=0$, the constant $C_n(\gamma)$ can be determined by
evaluating (\ref{algorithm:psn}) and the Legendre polynomial of degree $n$ at $0$.  We note that well-known
formulas (such as those found on page 145 of \cite{HTFI}) make it easy to evaluate $P_n(0)$ in time independent
of $n$.  Once $C_n(\gamma)$ has been determined,
the value of $\PS{n}(z;\gamma)$ can be calculated at any $z$ in the interval $[0,1-\exp(-\beta))$ by first
using the piecewise Chebyshev expansion (\ref{algorithm:chebexp}) to evaluate 
\begin{equation}
\PsiP{n}(x;\gamma) \ \ \mbox{and}\ \ \  \frac{d\PsiP{n}}{dx}(x;\gamma)
\end{equation}
and then plugging the resulting values into (\ref{algorithm:psn}).  
When $z$ is in  $\left(-1-\exp(-\beta),0\right)$, we use the well-known symmetry properties  of
$\PS{n}(z;\gamma)$ (it is an even function when $n$ is even and an odd function when $n$ is odd)
together with the procedure just described to evaluate it.

In the remainder of this section, we describe the method we use to solve ordinary differential
equations, the technique used to construct the expansion (\ref{algorithm:chebexp})
and an accelerated version of the algorithm obtained using the results of \cite{prolates2}.

\begin{subsection}{Adaptive solution of ordinary differential equations}

The algorithm of this paper entails solving several ordinary differential equations.
We use a fairly standard adaptive Chebyshev spectral solver to do so.  We now briefly describe its operation
in the case of the initial value problem
\begin{equation}
\left\{
\begin{aligned}
\bm{y}'(t) &= F(t,\bm{y}(t)), \ \ \ a < t < b,\\
\bm{y}(a) &= \bm{v}
\end{aligned}
\right.
\label{algorithm:system}
\end{equation}
where $F:\mathbb{R}^{n+1} \to \mathbb{R}^n$ is smooth and $\bm{v} \in \mathbb{R}^n$.
The solver can be easily modified to apply to a terminal value problem.

The solver takes as input a positive integer $k$, an interval $(a,b)$, a subroutine for evaluating the function $F$
and the vector $\bm{v}$.  It outputs $n$ piecewise $k^{th}$ order Chebyshev expansions,
one for each of the components $y_i(t)$ of the solution $\bm{y}$ of (\ref{algorithm:system}).

The solver maintains two lists of subintervals of $(a,b)$: one consisting of accepted intervals
and the other of intervals which have yet to processed.
Initially, the list of accepted intervals is empty and the list of 
intervals to process contains the single interval $(a,b)$
It then operates as follows until the list of intervals to process is empty:
\begin{enumerate}

\item
Find, in the list of interval to process, the interval $(c,d)$ such that
$c$ is as small as possible and remove this interval from the list.

\item
Solve the initial value problem
\begin{equation}
\left\{
\begin{aligned}
\bm{u}'(t) &= F(t,\bm{u}(t)), \ \ \ c< t < d,\\
\bm{u}(c) &= \bm{w}
\end{aligned}
\right.
\label{algorithm:ivp2}
\end{equation}
If $(c,d) = (a,b)$, then we take $\bm{w}=\bm{v}$.  Otherwise,
the value of the solution at the point $c$ has already been approximated, and we use that estimate
for $\bm{w}$ in (\ref{algorithm:ivp2}).

If the problem is linear, a straightforward Chebyshev spectral method (see, for instance, \cite{Trefethen})
is used to solve (\ref{algorithm:ivp2}).  Otherwise, 
the trapezoidal method (see, for instance, \cite{Ascher}) is first used to produce an initial
approximation $y_0$ of the solution and then Newton's method is applied to refine it.
The linearized problems are solved using a straightforward Chebyshev spectral method.

In any event, the result is a set of $k^{th}$ order Chebyshev expansions 
\begin{equation}
u_i(x)  \approx \sum_{j=0}^k \lambda_{ij}\ T_j\left(\frac{2}{d-c} x + \frac{c+d}{c-d}\right),\ \ \ i=1,\ldots,n,
\label{algorithm:exps}
\end{equation}
approximating  the components $u_1,\ldots,u_n$ of the solution of (\ref{algorithm:ivp2}).

\item
Compute the quantities
\begin{equation}
\frac{\sqrt{\sum_{j=k/2}^k \lambda_{ij}^2}}{\sqrt{\sum_{j=0}^k \lambda_{ij}^2}}, \ \ \ i=1,\ldots,n,
\end{equation}
where the $\lambda_{ij}$ are the coefficients in the expansions (\ref{algorithm:exps}).
If any of the resulting values is  larger than $100 \times \epsilon_0$, where $\epsilon_0$ denotes machine zero for 
the IEEE double precision number scheme,
then we split the interval into two halves $\left(c,\frac{c+d}{2}\right)$ and 
$\left(\frac{c+d}{2},d\right)$ and place them on the list of intervals to process.  Otherwise, we place the interval
$(c,d)$ on the list of accepted intervals.

\end{enumerate} 

At the conclusion of this procedure,  we have $k^{th}$ order Chebyshev expansions
for each component of the solution, with the list of accepted intervals determining the
partition for each expansion.
\end{subsection}


\begin{subsection}{Construction of $\PsiP{n}(x;\gamma)$}

To construct $\PsiP{n}(x;\gamma)$, we first calculate the values of the function $M_\nu(x;\gamma)$ 
and its first two derivatives with respect to $x$ at $0$. We do this by evaluating the logarithmic derivative 
\begin{equation}
s(t) = \frac{f'(t)}{f(t)}
\end{equation}
of the function
\begin{equation}
f(t) = S^{(3)}_n(it;\gamma)\sqrt{1+t^2}
\end{equation}
and its derivative at the point $0$.    Since $f$ is a solution of (\ref{spheroidal:irswe}),  $s$ satisfies
\begin{equation}
s'(t) + (s(t))^2 + q_1(t) = 0,
\end{equation}
where $q_1$ is the coefficient in (\ref{spheroidal:irswe}) with $\chi$ taken to be equal to the Sturm-Liouville eigenvalue
$\chi_n(\gamma)$.    Moreover, from (\ref{spheroidal:s3asym}),  it is apparent that
\begin{equation}
s(t) = -\gamma + \mathcal{O}\left(\frac{1}{t}\right).
\end{equation}
Accordingly, we construct $s$ by adaptively solving the terminal value problem
\begin{equation}
\left\{
\begin{aligned}
s'(t) + (s(t))^2 + q_1(t) &= 0,\ \ \ 0 < t < c\\
s(c) &= -\gamma,
\end{aligned}
\right.
\label{algorithm:axis}
\end{equation}
where $c$ is a suitable large constant (we take $c=10^{30}$ in the experiments presented here).  Because of Conjecture~\ref{conjecture2},
we expect the function $s$ to be well-behaved on the interval $(0,\infty)$ and for the cost of solving
(\ref{algorithm:axis}) to be slowly growing with $\gamma$ and $n$.

Using the relations found in Section~\ref{section:riccati} and the definitions of Section~\ref{section:spheroidal},
the values of  $w(x) = M_\nu(1-\exp(-x);\gamma)(2-\exp(-x))$ and its first two derivatives at $x=0$ can
be expressed in terms of $s(0)$ and $s'(0)$.  Indeed, $w$ is the solution of the initial value problem
\begin{equation}
\left\{
\begin{aligned}
&w'''(x) + 4 q_2(x) w'(x) + 2q_2'(x) w(x) = 0, \ \ \  0 \leq x < b, \\
&w(0) =  -\frac{1}{s(0)}\\
&w'(0) = -\frac{1}{s(0)}\\
&w''(0) = -\frac{1}{s(0)}+2\frac{s'(0)}{s(0)},
\end{aligned}
\right.
\label{algorithm:ivp}
\end{equation}
where $q_2$ is the coefficient in (\ref{spheroidal:erswe}) with $\chi$ taken to be equal to $\chi_n(\gamma)$.   
The next step of our algorithm
consists of  solving (\ref{algorithm:ivp})
to construct a piecewise $k^{th}$ order Chebyshev expansion for $M_n(1-\exp(-x);\gamma)(2-\exp(-x))$.
We let $0 = x_0 < x_1 < x_2 < \ldots < x_n=b$ be the partition on which
this piecewise expansion is given.  

Since
%
\begin{equation}
\begin{aligned}
\PsiP{n}(x;\gamma) &= \int_{\infty}^x \frac{\gamma}{M_n(1-\exp(-u);\gamma)(2-\exp(-u))}\ du \\
&\approx \int_{\beta}^x \frac{\gamma}{M_n(1-\exp(-u);\gamma)(2-\exp(-u))}\ du,
\end{aligned}
\label{algorithm:intleft}
\end{equation}
spectral integration can be used to construct a $k^{th}$ order  piecewise Chebyshev expansion 
for $\PsiP{n}(x;\gamma)$.  More explicitly, the final step of our algorithm consists
of  traversing the intervals $(x_{i-1},x_{i})$ in decreasing order (i.e.,
$i=n,n-1,n\ldots,1$) and, on each interval $(x_{i-1},x_i)$, constructing a $k^{th}$ order
Chebyshev expansion representing $\PsiP{n}(x;\gamma)$  
by applying a $k\times k$ spectral integration
matrix to the vector consisting of the coefficients in the Chebyshev expansion of 
 $M_n(1-\exp(-x);\gamma)(2-\exp(-x))$ over the interval $(x_{i-1},x_i)$.  The result is the
piecewise Chebyshev expansion (\ref{algorithm:chebexp}) for $\PsiP{n}(x;\gamma)$.

The derivative of the phase function decays extremely slowly to $0$ 
when the parameters $\gamma$ and $n$ are of small magnitude,  and, in this event, it is necessary to choose $\beta$ to be extremely close
to $1$ to achieve high accuracy in (\ref{algorithm:intleft}).  This is  what motivates the exponential change of variables used
to obtain the form (\ref{spheroidal:erswe}) of the reduced spheroidal wave equation and the decision to make  $\beta$ 
so large.  We note that for most values of the parameters, $\beta$ can be taken to be much smaller without losing accuracy.

We prefer to solve Appell's equation (\ref{riccati:appell})  over Kummer's equation
(\ref{riccati:kummer}) because of difficulties which are encountered when solving
the latter numerically.  The solution of Kummer's equation is the derivative
of the phase function $\PsiP{n}(x;\gamma)$.  In cases in which (\ref{spheroidal:erswe}) has
a turning point, the derivative of $\PsiP{n}(x;\gamma)$ decays rapidly to $0$ once
the turning point is reached.  This creates numerical complications when solving
Kummer's equation owing to the presence of the derivative of the phase function
in the denominators of the some of the terms in (\ref{riccati:kummer}).
On the other hand, the modulus function which satisfies Appell's equation
increasing on  $(0,\beta)$ (indeed, it is absolutely monotone if our conjectures are correct).

\end{subsection}

\vskip 1em

\begin{subsection}{Accelerated version of the algorithm}

In \cite{prolates2}, a numerical method for evaluating $\chi_n(\gamma)$ in time 
which is bounded independent of $n$ and $\gamma$ is described.  It also allows
for the $\mathcal{O}(1)$ calculation of the values of 
$\PsiP{n}(z;\gamma)$  and its first few derivatives at the point $z=0$.
Using these quantities, the 
values of modulus function $M_n(1-\exp(-x);\gamma)(2-\exp(-x))$ 
and its first two derivatives at $x=0$ can be easily obtained.
The algorithm of this paper can be significantly accelerated by exploiting this capability.
First, when the values of the modulus function   $M_n(1-\exp(-x);\gamma)(2-\exp(-x))$ and its first two 
derivatives at $x=0$ are known, there is no  need to solve  the terminal value problem (\ref{algorithm:axis}).
Second, when the value of $\PsiP{n}(0;\gamma)$ is known, Formula~(\ref{algorithm:intleft})
can be replaced with 
\begin{equation}
\begin{aligned}
\PsiP{n}(x;\gamma) &= \PsiP{n}(0;\gamma) + \int_{0}^x \frac{\gamma}{M_n(1-\exp(-u);\gamma)(2-\exp(-u))}\ du.
\end{aligned}
\label{algorithm:intright}
\end{equation}
Since we are no longer using (\ref{spheroidal:psipconstant}) to determine the constant
of integration, it is no longer necessary to 
calculate $\PsiP{n}(x;\gamma)$ on a large interval.  We instead construct a piecewise
$k^{th}$ Chebyshev expansion for it on the interval $[0,30)$.  This allows for the evaluation
of $\PS{n}(z;\gamma)$ for $z$ in $[0,1-\exp(-30))$, which more than suffices for most purposes.

\end{subsection}

\end{section}

\begin{figure}[!t]
\hfil
\includegraphics[width=.49\textwidth]{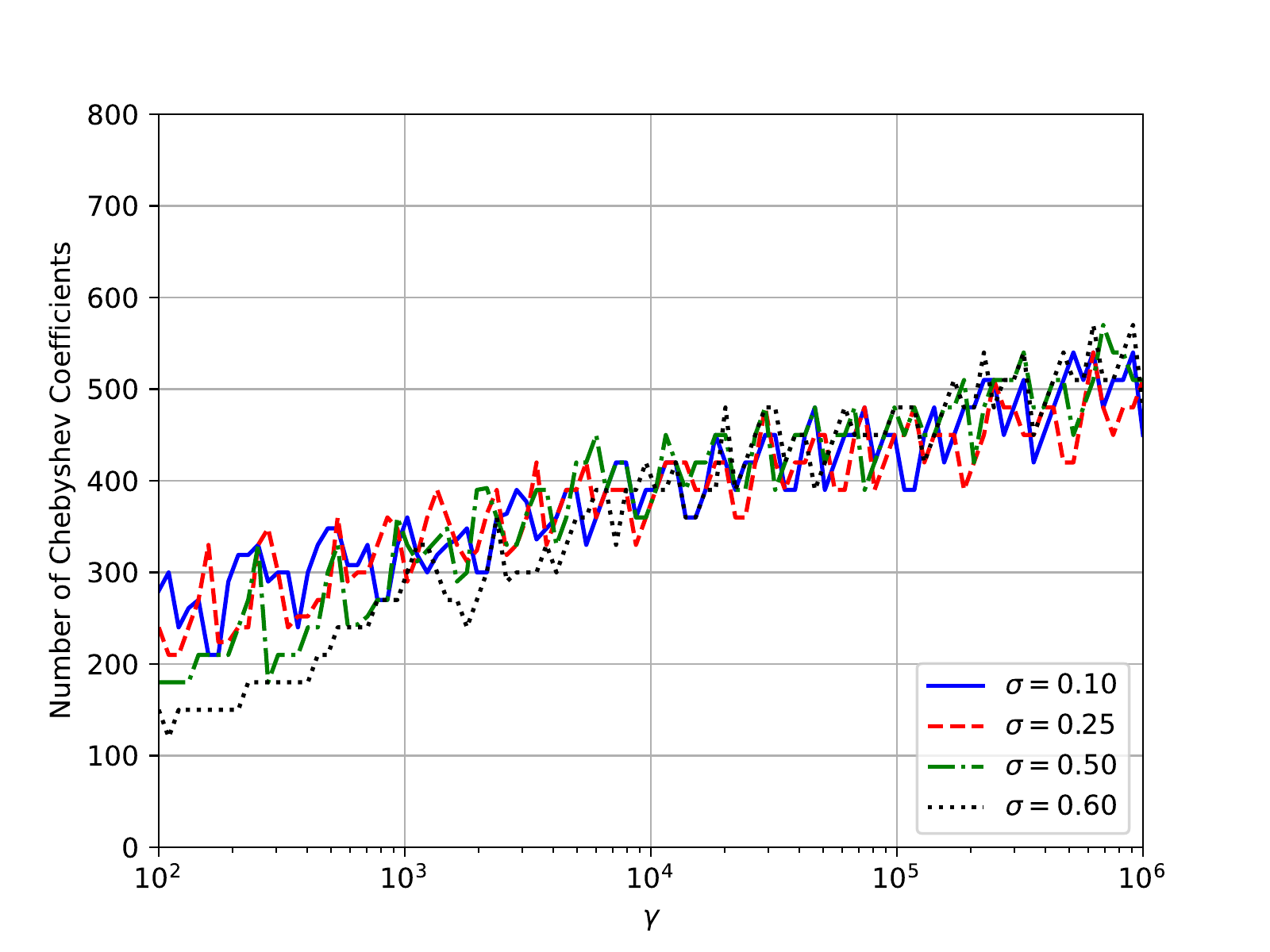}
\hfil
\includegraphics[width=.49\textwidth]{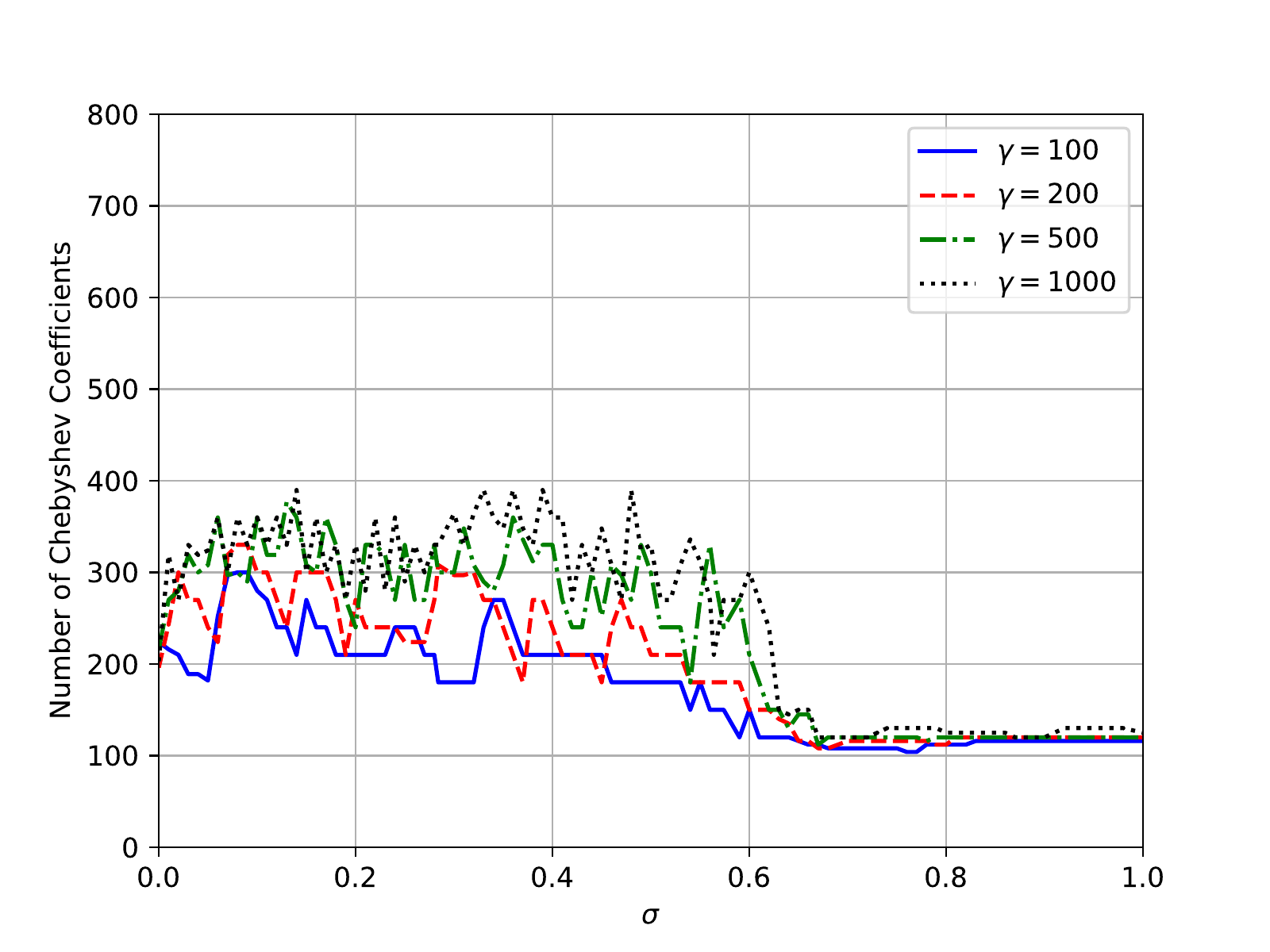}
\hfil

\vskip 3em

\hfil
\includegraphics[width=.49\textwidth]{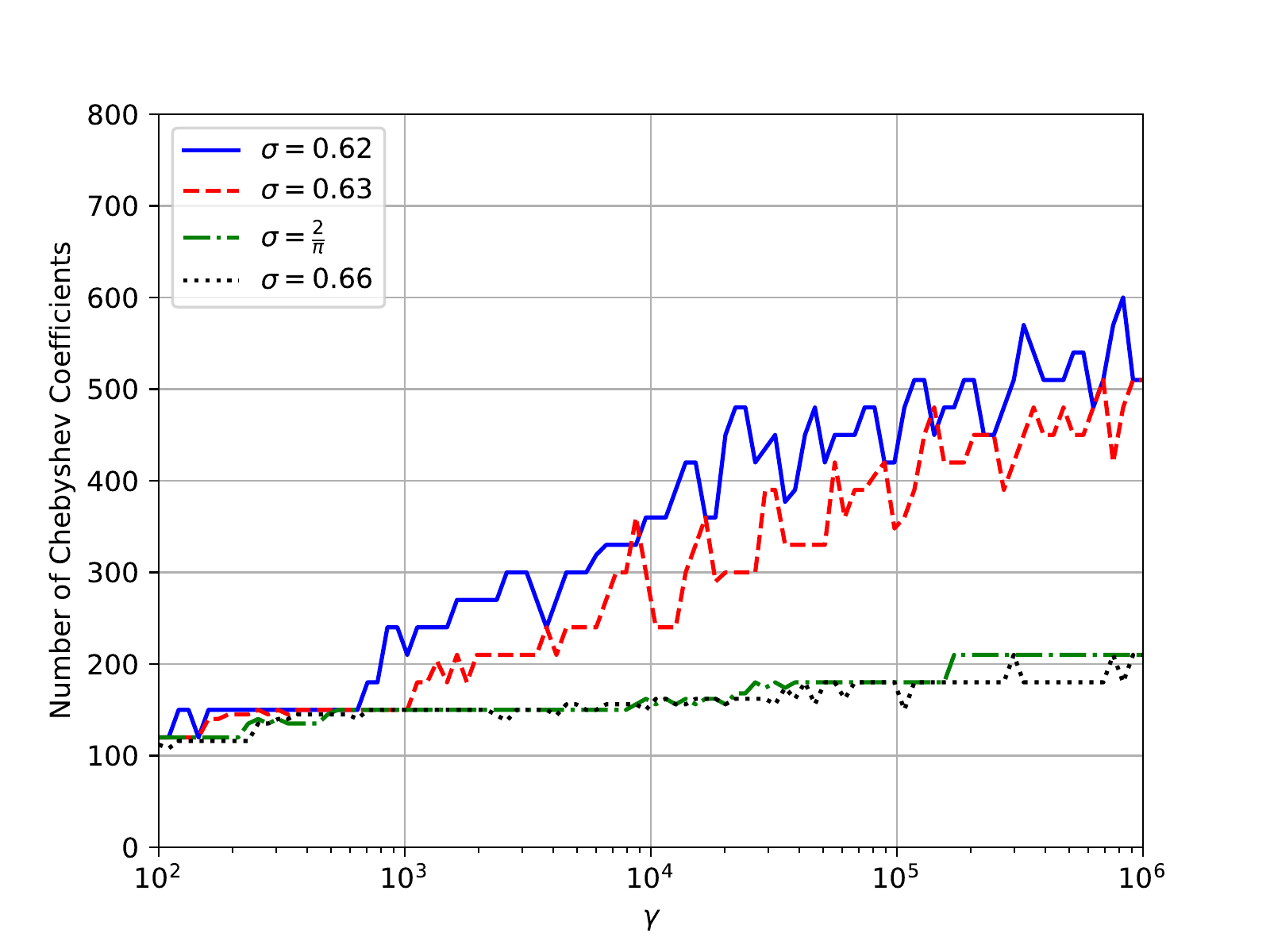}
\hfil
\includegraphics[width=.49\textwidth]{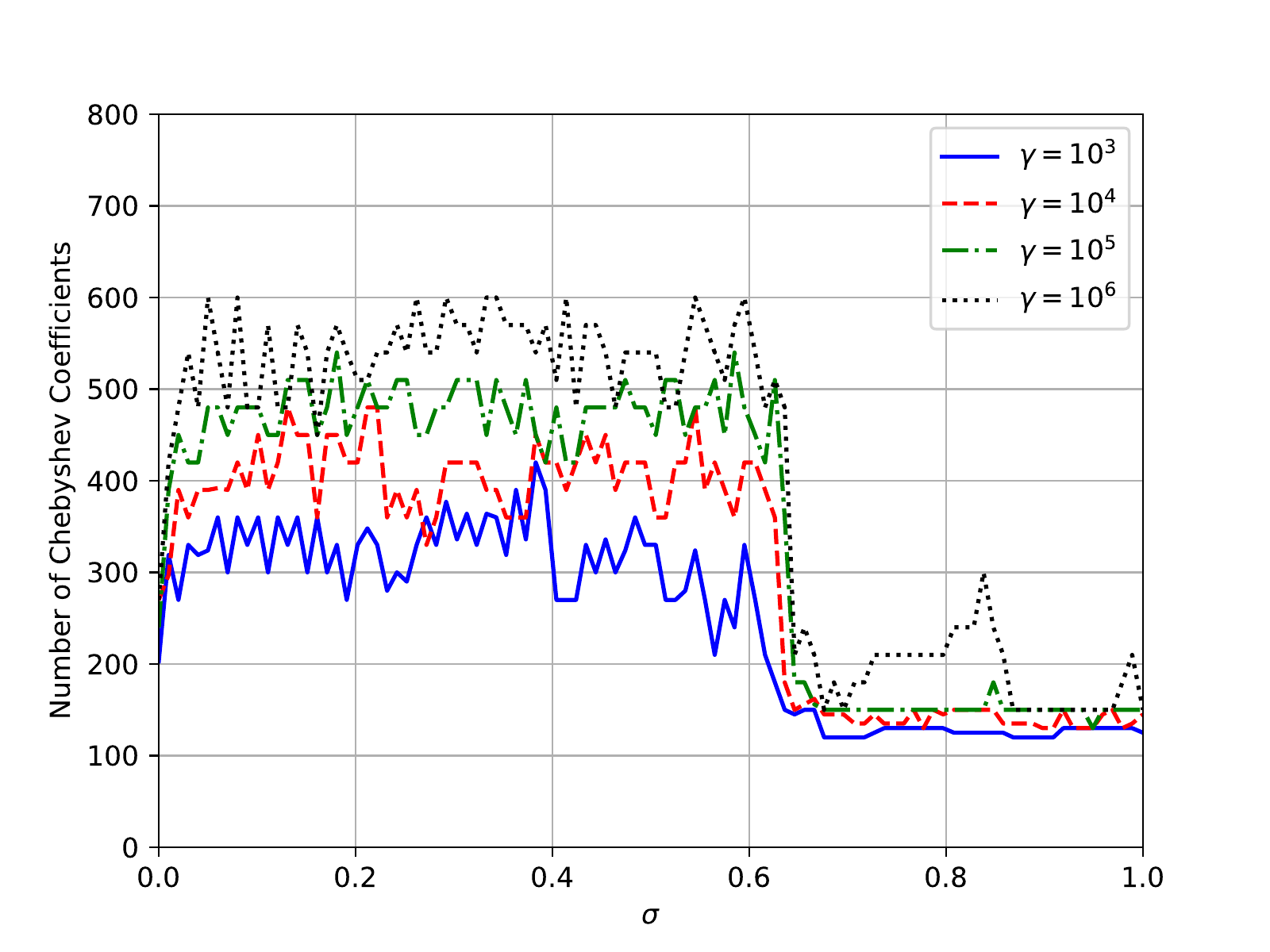}
\hfil

\vskip 3em

\hfil
\includegraphics[width=.49\textwidth]{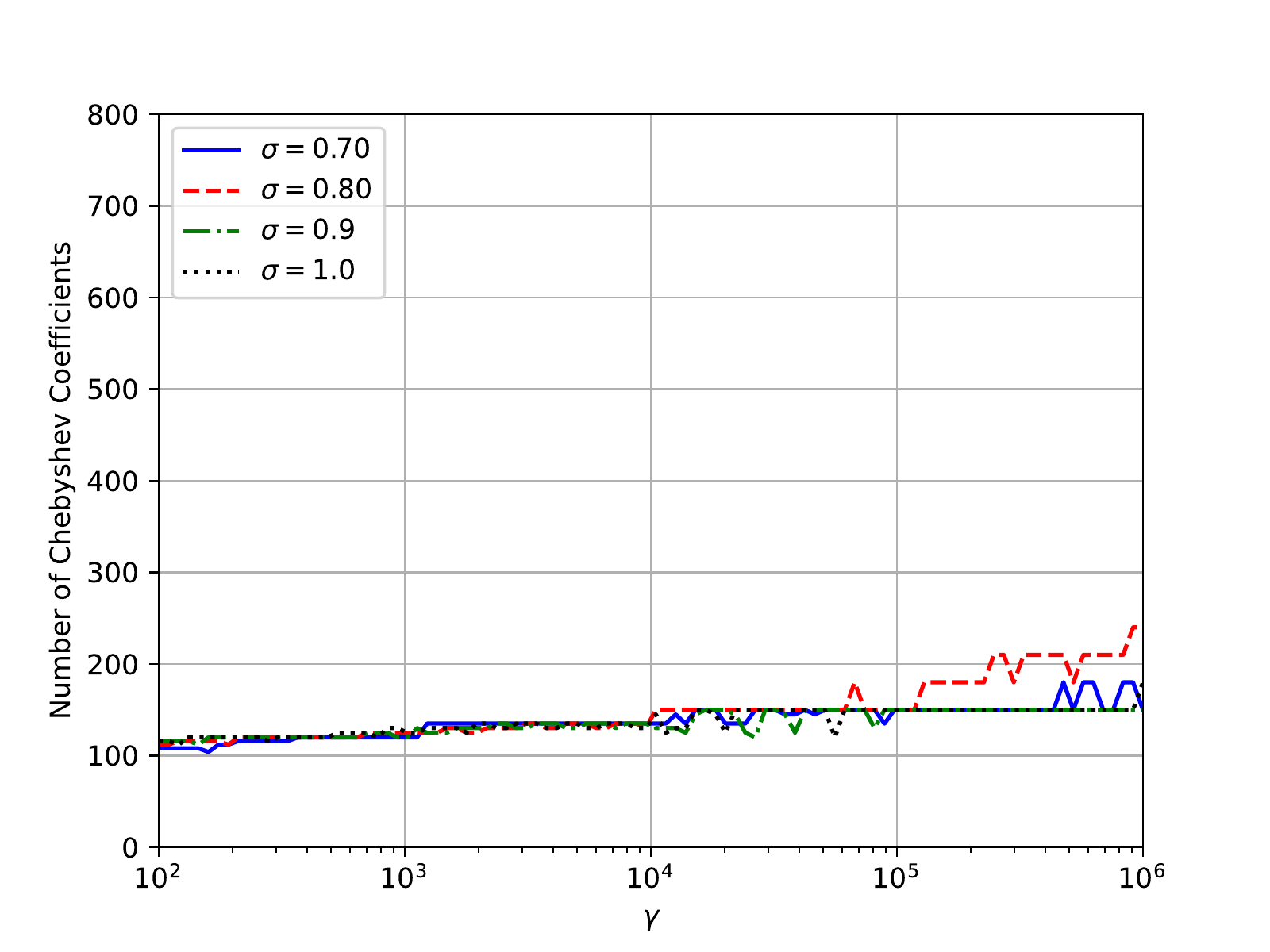}
\hfil
\includegraphics[width=.49\textwidth]{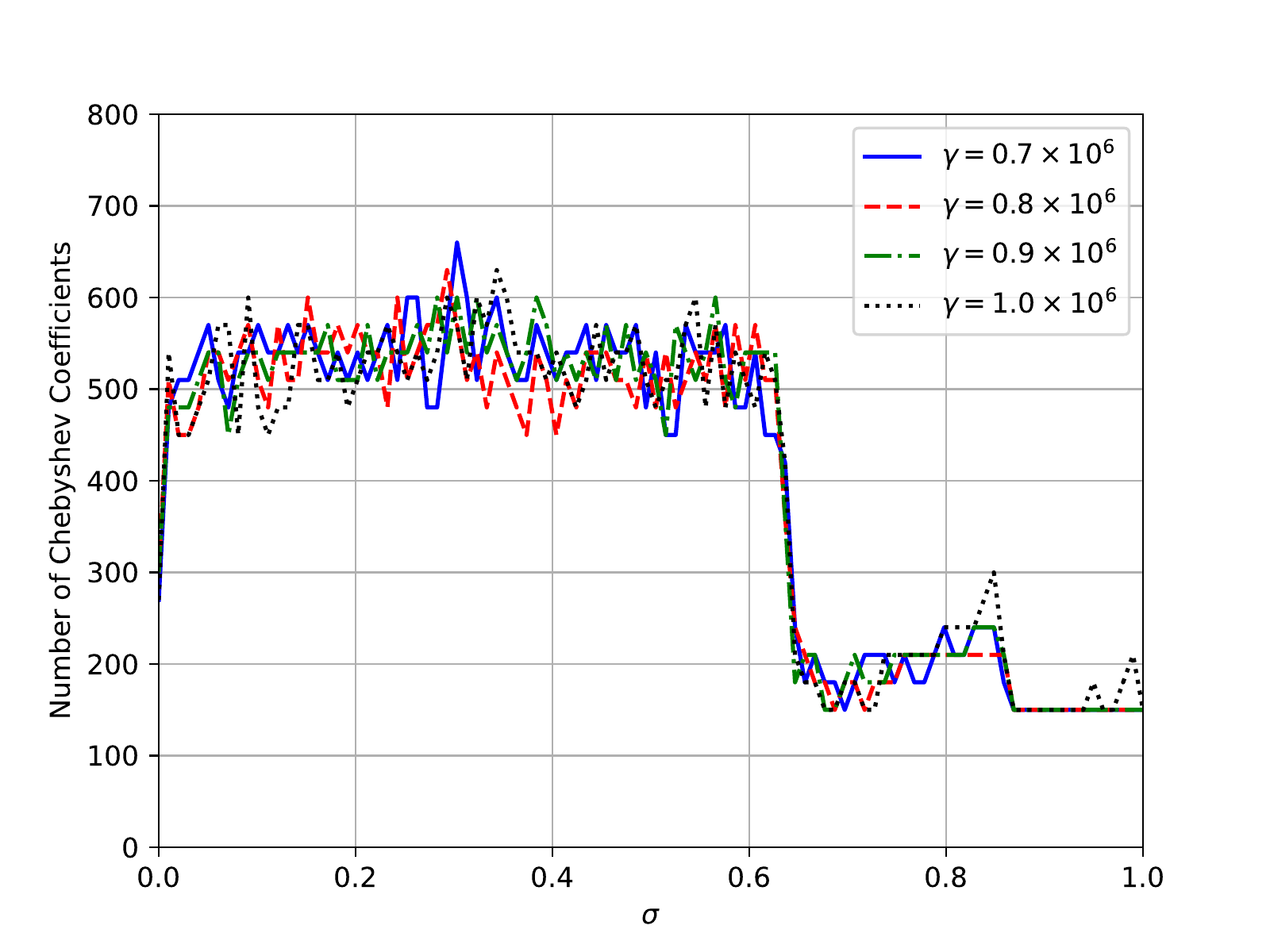}
\hfil

\caption{The number of Chebyshev coefficients required to represent $\PsiP{n}(\gamma)$.
Each of the graphs on the left gives the number of coefficients as function of $\gamma$ for several
values of $\sigma$, while the plots on the right give the number of coefficients
as a function of $\sigma$ for several values of $\gamma$.  
A logarithmic scale is used for the x-axis in each of the plots on the left.}
\label{experiments:figure1}
\end{figure}

\begin{section}{Numerical Experiments}

In this section, we present the results of numerical experiments which were conducted
to illustrate the effectiveness of the algorithm of this article and to 
measure the dependence of the cost of representing $\PsiP{n}(x;\gamma)$ 
on the parameters.  As discussed in Section~\ref{section:spheroidal}, for the  purposes of numerical
computation, the effective dimension of the space of functions with bandlimit $\gamma$ is 
\begin{equation}
\frac{2}{\pi} \gamma + \mathcal{O}\left(\log(\gamma)\right);
\end{equation}
in particular, the number of prolate functions of interest depends on the bandlimit $\gamma$.  
Moreover, the qualitative behaviour of $\PS{n}(z;\gamma)$ is best understood in terms of the
ratio of $n$ to $\gamma$.  Accordingly, we introduce a new parameter
$\sigma$ which is related to $n$ via
\begin{equation}
n = \mbox{round} \left(\gamma \sigma\right),
\end{equation}
and the results presented here are discussed in terms of the parameters
$\gamma$ and $\sigma$ rather than $\gamma$ and $n$.

The code for our experiments was written in Fortran and compiled with version 11.1.0 of the GNU
Fortran compiler.  They were performed on a desktop computer equipped with an AMD Ryzen 3900x
processor and 32GB of RAM.  An implementation of our algorithm and  code for 
conducting all of the experiments discussed here is available on GitHub at the following address:
\begin{center}
\url{https://github.com/JamesCBremerJr/Prolates}
\end{center}
%

In some of these experiments, we compared the performance of our algorithm
with that of the Osipov-Xiao-Rokhlin method \cite{Osipov-Rokhlin-Xiao}.  Its running time is highly
dependent on the dimension of the tridiagonal matrix formed in order to calculate
the coefficients in the Legendre expansion of $\PS{n}(z;\gamma)$.  Most implementations 
use a highly conservative value for this dimension.  The authors of \cite{Osipov-Rokhlin-Xiao}, for instance,
take it to be  $1000 + n + \left \lfloor 1.1 \gamma \right\rfloor$ in their implementation.
In many cases, $\PS{n}(z;\gamma)$ can be represented much more efficiently than this.
The experiments of  \cite{Feichtinger}, though, suggest that the necessary
dimension grows as  $\mathcal{O}\left(n + \sqrt{n\gamma}\right)$.
It is difficult, however, to find a simple formula which suffices in all cases
of interest.  Accordingly,  our implementation of the Osipov-Rokhlin-Xiao algorithm initially
takes the dimension to be 
\begin{equation}
50 + \left\lfloor\frac{2}{\pi}n\right\rfloor +  \left\lfloor \sqrt{\gamma n}\right\rfloor,
\end{equation}
which we found to be sufficient for a large range of parameters, and 
then increases it adaptively as needed to ensure high accuracy.
Our implementation can be found in the GitHub repository mentioned above.


\begin{subsection}{The cost of representing $\PsiP{n}(x;\gamma)$}

In the experiments discussed here, we measured the  cost to represent $\PsiP{n}(x;\gamma)$ 
on the interval $[0,30)$
using a piecewise $k^{th}$ order Chebyshev expansion with $k=29$.    Table~\ref{experiments:table1} and 
Figure~\ref{experiments:figure1} show the results. 

\begin{table}[!h]
\begin{center}
\begin{tabular}{ccr@{\hspace{2em}}ccr}
\toprule
Range of $\gamma$ & Range of $\sigma$ & Max Coefs  &Range of $\gamma$ & Range of $\sigma$ & Max Coefs \\\midrule
100 to 500                                                                                           &$0.00 - 0.25 $ & $  420$ &10,000 to 50,000                                                                                     &$0.00 - 0.25 $ & $  540$ \\
\addlinespace[.125em]
 & 
$0.25 - 0.50 $ & $  420$ & & 
$0.25 - 0.50 $ & $  570$ \\
\addlinespace[.125em]
 & 
$0.50 - 0.75 $ & $  360$ & & 
$0.50 - 0.75 $ & $  570$ \\
\addlinespace[.125em]
 & 
$0.75 - 1.00 $ & $  120$ & & 
$0.75 - 1.00 $ & $  180$ \\
\addlinespace[.125em]
\addlinespace[.25em]
500 to 1,000                                                                                         &$0.00 - 0.25 $ & $  420$ &50,000 to 100,000                                                                                    &$0.00 - 0.25 $ & $  570$ \\
\addlinespace[.125em]
 & 
$0.25 - 0.50 $ & $  450$ & & 
$0.25 - 0.50 $ & $  570$ \\
\addlinespace[.125em]
 & 
$0.50 - 0.75 $ & $  364$ & & 
$0.50 - 0.75 $ & $  570$ \\
\addlinespace[.125em]
 & 
$0.75 - 1.00 $ & $  130$ & & 
$0.75 - 1.00 $ & $  180$ \\
\addlinespace[.125em]
\addlinespace[.25em]
1,000 to 5,000                                                                                       &$0.00 - 0.25 $ & $  450$ &100,000 to 500,000                                                                                   &$0.00 - 0.25 $ & $  630$ \\
\addlinespace[.125em]
 & 
$0.25 - 0.50 $ & $  480$ & & 
$0.25 - 0.50 $ & $  630$ \\
\addlinespace[.125em]
 & 
$0.50 - 0.75 $ & $  450$ & & 
$0.50 - 0.75 $ & $  630$ \\
\addlinespace[.125em]
 & 
$0.75 - 1.00 $ & $  145$ & & 
$0.75 - 1.00 $ & $  240$ \\
\addlinespace[.125em]
\addlinespace[.25em]
5,000 to 10,000                                                                                      &$0.00 - 0.25 $ & $  510$ &500,000 to 1,000,000                                                                                 &$0.00 - 0.25 $ & $  630$ \\
\addlinespace[.125em]
 & 
$0.25 - 0.50 $ & $  510$ & & 
$0.25 - 0.50 $ & $  660$ \\
\addlinespace[.125em]
 & 
$0.50 - 0.75 $ & $  480$ & & 
$0.50 - 0.75 $ & $  630$ \\
\addlinespace[.125em]
 & 
$0.75 - 1.00 $ & $  150$ & & 
$0.75 - 1.00 $ & $  300$ \\
\addlinespace[.125em]
\addlinespace[.25em]
\bottomrule
\end{tabular}

\end{center}
\caption{The cost of representing $\PsiP{n}(x;\gamma)$ for various ranges of values of $\gamma$
and $\sigma$.  
To generate the data presented in each row of this table, $100$ equispaced
values of $\gamma$ and $\sigma$ in the ranges indicated were sampled and a
representation of $\PsiP{n}(x;\gamma)$ over the interval $[0,1-\exp(-30))$
was constructed for each of the $10,000$ resulting pairs of the parameters.}
\label{experiments:table1}
\end{table}

 Table~\ref{experiments:table1} reports the number of 
coefficients required to represent the phase function on the interval $[0,1-\exp(-30))$ for 
various ranges of values of $\gamma$ and $\sigma$.  For each row of the table,
we first sampled $100$ equispaced values of $\gamma$
and $100$ equispaced values of $\sigma$ in the indicated ranges.
Then, we constructed a piecewise Chebyshev  expansion of $\PsiP{n}(x;\gamma)$ for each of the $10,000$ pairs of sampled
values and determined the number of coefficients in the largest expansion encountered.

Each plot on the left side  of Figure~\ref{experiments:figure1} 
gives the number of Chebyshev coefficients in the expansion representing $\PsiP{n}(x;\gamma)$ as  
function of $\gamma$ for various values of $\sigma$, while
each plot on the right side shows the number of Chebyshev coefficients in the expansion
as a function of $\sigma$ for various values of $\gamma$.  A logarithmic scale is used for
the x-axis in each plot on the left.

We immediately draw several conclusions from these experiments.  First, we note that 
no more than $800$ Chebyshev coefficients were required to represent any of the expansions
formed during the course of these experiments.  
Second, we observe that while
the cost to represent $\PsiP{n}(x;\gamma)$ grows with $\gamma$, it does so at a modest rate.
Indeed, the plots on the left side of  Figure~\ref{experiments:figure1} indicate that
for values of $\sigma$ somewhat larger than $2/\pi$, the number of coefficients is essentially
constant as a function of $\gamma$, while it increases sublogarithmically with
$\gamma$ for $\sigma$ which are less than $2/\pi$. 
Finally, we observe from the plots on the right side of Figure~\ref{experiments:figure1}
that a rapid drop in the cost of representing the phase function
occurs when $\sigma$ is close to $2/\pi$.  
 For small values of $\sigma$, the reduced spheroidal
wave equation has turning points on $(0,1)$ and its solutions are oscillatory
only on part of that interval.  However, starting when $\sigma$ is a bit larger than $2/\pi$,
the solutions of the reduced spheroidal wave equation are oscillatory on all of $(0,1)$.
Evidently, the cost of representing $\PsiP{n} (x;\gamma)$ in the oscillatory
regime  is bounded independent of $\gamma$, while the cost of representing
it in the nonoscillatory regime grows sublogarithmically with $\gamma$.

\end{subsection}

%
%
\begin{subsection}{The accuracy with which $\PS{n}(z;\gamma)$ is evaluated}

We now describe experiments conducted to measure the  accuracy with which the algorithm of this paper
evaluates the angular spheroidal wave functions of the first kind of order zero by comparison  with 
the Osipov-Xiao-Rokhlin method.  Table~\ref{experiments:table2} and Figure~\ref{experiments:figure2}
give the results.

Table~\ref{experiments:table2} gives the accuracy  of our algorithm for 
certain ranges of the parameters $\gamma$ and $\sigma$.
To generate the data presented in each row of this table, we first sampled $100$ equispaced
values of $\gamma$ and $\sigma$ in the ranges indicated. Then, for each of the $10,000$ pairs of sampled values,   we 
evaluated $\PS{n}(z;\gamma)$ at $100$ equispaced
points $\{x_k\}$ in the interval $(0,1)$ using the algorithm of this paper and the Osipov-Xiao-Rokhlin
method.  The largest absolute error observed is reported in Table~\ref{experiments:table2}.

Figure~\ref{experiments:figure2} contains plots showing the dependence of the accuracy of our algorithm
on $\gamma$ and $\sigma$.  For each pair of parameters considered, 
$\PS{n}(z;\gamma)$ was evaluated at $100$ equispaced points $\{x_k\}$ on the interval $(0,1)$ 
and the largest absolute error was determined.
Each plot on the left side  of that figure  gives the accuracy with which $\PS{n}(z;\gamma)$ is evaluated as 
function of $\gamma$ for various values of $\sigma$, while
each plot on the right side plots the accuracy 
as a function of $\sigma$ for various values of $\gamma$. 
A logarithmic scale is used for the x-axis in each plot on the left.

As expected, accuracy is lost as both $\gamma$ and $\sigma$ increase.  This occurs because
the magnitude of $\PsiP{n}(x;\gamma)$ increases with both of these variables and the accuracy
with which the sine function is evaluated decreases with the magnitude of its argument.
We note that the condition number of the reduced spheroidal wave equation and the condition
number of evaluation of its solutions increases with the parameters $\gamma$ and $\sigma$
as well, so some loss of precision is expected.  We also observe that the largest error observed
during the course of these experiments was approximately $10^{-10}$, and this was
only for values of $\gamma$ near $1,000,000$.

\begin{table}[!t]
\begin{center}
\begin{tabular}{ccc@{\hspace{1em}}ccc}
\toprule
Range of $\gamma$ & Range of $\sigma$ & Max Error  &Range of $\gamma$ & Range of $\sigma$ & Max Error \\\midrule
100 to 500                                                                                           &$0.00 - 0.25 $ & $6.74\e{-14}$ &10,000 to 50,000                                                                                     &$0.00 - 0.25 $ & $5.90\e{-13}$ \\\addlinespace[.125em]
 & 
$0.25 - 0.50 $ & $1.02\e{-13}$ & & 
$0.25 - 0.50 $ & $9.07\e{-13}$ \\\addlinespace[.125em]
 & 
$0.50 - 0.75 $ & $3.91\e{-13}$ & & 
$0.50 - 0.75 $ & $4.88\e{-12}$ \\\addlinespace[.125em]
 & 
$0.75 - 1.00 $ & $1.45\e{-13}$ & & 
$0.75 - 1.00 $ & $2.10\e{-12}$ \\\addlinespace[.125em]
\addlinespace[.25em]
500 to 1,000                                                                                         &$0.00 - 0.25 $ & $1.02\e{-13}$ &50,000 to 100,000                                                                                    &$0.00 - 0.25 $ & $8.62\e{-13}$ \\\addlinespace[.125em]
 & 
$0.25 - 0.50 $ & $1.49\e{-13}$ & & 
$0.25 - 0.50 $ & $1.42\e{-12}$ \\\addlinespace[.125em]
 & 
$0.50 - 0.75 $ & $7.02\e{-13}$ & & 
$0.50 - 0.75 $ & $2.12\e{-11}$ \\\addlinespace[.125em]
 & 
$0.75 - 1.00 $ & $2.09\e{-13}$ & & 
$0.75 - 1.00 $ & $2.03\e{-11}$ \\\addlinespace[.125em]
\addlinespace[.25em]
1,000 to 5,000                                                                                       &$0.00 - 0.25 $ & $2.09\e{-13}$ &100,000 to 500,000                                                                                   &$0.00 - 0.25 $ & $1.94\e{-12}$ \\\addlinespace[.125em]
 & 
$0.25 - 0.50 $ & $3.25\e{-13}$ & & 
$0.25 - 0.50 $ & $3.39\e{-12}$ \\\addlinespace[.125em]
 & 
$0.50 - 0.75 $ & $2.14\e{-12}$ & & 
$0.50 - 0.75 $ & $5.32\e{-11}$ \\\addlinespace[.125em]
 & 
$0.75 - 1.00 $ & $5.54\e{-13}$ & & 
$0.75 - 1.00 $ & $5.44\e{-11}$ \\\addlinespace[.125em]
\addlinespace[.25em]
5,000 to 10,000                                                                                      &$0.00 - 0.25 $ & $2.76\e{-13}$ &500,000 to 1,000,000                                                                                 &$0.00 - 0.25 $ & $2.87\e{-12}$ \\\addlinespace[.125em]
 & 
$0.25 - 0.50 $ & $4.31\e{-13}$ & & 
$0.25 - 0.50 $ & $4.51\e{-12}$ \\\addlinespace[.125em]
 & 
$0.50 - 0.75 $ & $3.13\e{-12}$ & & 
$0.50 - 0.75 $ & $1.09\e{-10}$ \\\addlinespace[.125em]
 & 
$0.75 - 1.00 $ & $1.10\e{-12}$ & & 
$0.75 - 1.00 $ & $7.52\e{-11}$ \\\addlinespace[.125em]
\addlinespace[.25em]
\bottomrule
\end{tabular}

\end{center}
\caption{The accuracy with which $\PS{n}(z;\gamma)$ is evaluated for various ranges of values of $\gamma$
and $\sigma$.  
}
\label{experiments:table2}
\end{table}

\end{subsection}

%
%
\begin{subsection}{The time required to construct $\PsiP{n}(x;\gamma)$}

In these experiments, we measured the time required to calculate the phase function
$\PsiP{n}(x;\gamma)$ using both the accelerated and unaccelerated algorithm described in this paper.  In some
of these experiments, we  compared it with the time taken by the Osipov-Xiao-Rokhlin method to construct 
the Legendre expansion representing  $\PS{n}(z;\gamma)$.
Table~\ref{experiments:table3} and Figures~\ref{experiments:figure3} and \ref{experiments:figure4}
give the results.

\begin{table}[!t]
\begin{center}
\begin{tabular}{ccccccr}
\toprule
Range of $\gamma$ & Range of $\sigma$ & Average Time   & Average time & Average time        & Ratio\\
              &                   & Unaccelerated  & Accelerated    & Rokhlin, et. al.\\
\midrule
100 to 500                                                                                           &$0.00 - 0.25 $ & $5.35\e{-03}$ &$3.05\e{-04}$ &$9.76\e{-05}$ &$    0.32
$ \\\addlinespace[.125em]
 & 
$0.25 - 0.50 $ & $3.32\e{-03}$ &$2.62\e{-04}$ &$1.31\e{-04}$ &$    0.50
$ \\\addlinespace[.125em]
 & 
$0.50 - 0.75 $ & $1.90\e{-03}$ &$1.16\e{-04}$ &$1.71\e{-04}$ &$    1.47
$ \\\addlinespace[.125em]
 & 
$0.75 - 1.00 $ & $1.85\e{-03}$ &$6.95\e{-05}$ &$2.11\e{-04}$ &$    3.04
$ \\\addlinespace[.125em]
\addlinespace[.25em]
500 to 1,000                                                                                         &$0.00 - 0.25 $ & $6.78\e{-03}$ &$3.54\e{-04}$ &$1.85\e{-04}$ &$    0.52
$ \\\addlinespace[.125em]
 & 
$0.25 - 0.50 $ & $5.42\e{-03}$ &$3.43\e{-04}$ &$2.87\e{-04}$ &$    0.84
$ \\\addlinespace[.125em]
 & 
$0.50 - 0.75 $ & $2.24\e{-03}$ &$1.53\e{-04}$ &$3.96\e{-04}$ &$    2.58
$ \\\addlinespace[.125em]
 & 
$0.75 - 1.00 $ & $1.97\e{-03}$ &$9.26\e{-05}$ &$4.90\e{-04}$ &$    5.29
$ \\\addlinespace[.125em]
\addlinespace[.25em]
1,000 to 5,000                                                                                       &$0.00 - 0.25 $ & $7.20\e{-03}$ &$4.04\e{-04}$ &$5.73\e{-04}$ &$    1.42
$ \\\addlinespace[.125em]
 & 
$0.25 - 0.50 $ & $6.94\e{-03}$ &$4.18\e{-04}$ &$1.09\e{-03}$ &$    2.61
$ \\\addlinespace[.125em]
 & 
$0.50 - 0.75 $ & $3.26\e{-03}$ &$2.11\e{-04}$ &$1.48\e{-03}$ &$    7.01
$ \\\addlinespace[.125em]
 & 
$0.75 - 1.00 $ & $2.16\e{-03}$ &$9.84\e{-05}$ &$1.87\e{-03}$ &$   19.08
$ \\\addlinespace[.125em]
\addlinespace[.25em]
5,000 to 10,000                                                                                      &$0.00 - 0.25 $ & $7.40\e{-03}$ &$4.36\e{-04}$ &$1.34\e{-03}$ &$    3.08
$ \\\addlinespace[.125em]
 & 
$0.25 - 0.50 $ & $7.27\e{-03}$ &$4.54\e{-04}$ &$2.67\e{-03}$ &$    5.88
$ \\\addlinespace[.125em]
 & 
$0.50 - 0.75 $ & $4.09\e{-03}$ &$2.47\e{-04}$ &$3.65\e{-03}$ &$   14.80
$ \\\addlinespace[.125em]
 & 
$0.75 - 1.00 $ & $2.47\e{-03}$ &$9.88\e{-05}$ &$4.63\e{-03}$ &$   46.87
$ \\\addlinespace[.125em]
\addlinespace[.25em]
10,000 to 50,000                                                                                     &$0.00 - 0.25 $ & $7.58\e{-03}$ &$4.74\e{-04}$ &$5.12\e{-03}$ &$   10.79
$ \\\addlinespace[.125em]
 & 
$0.25 - 0.50 $ & $7.53\e{-03}$ &$4.99\e{-04}$ &$1.04\e{-02}$ &$   20.96
$ \\\addlinespace[.125em]
 & 
$0.50 - 0.75 $ & $4.94\e{-03}$ &$2.94\e{-04}$ &$1.45\e{-02}$ &$   49.39
$ \\\addlinespace[.125em]
 & 
$0.75 - 1.00 $ & $3.12\e{-03}$ &$1.02\e{-04}$ &$1.81\e{-02}$ &$  177.52
$ \\\addlinespace[.125em]
\addlinespace[.25em]
50,000 to 100,000                                                                                    &$0.00 - 0.25 $ & $7.71\e{-03}$ &$5.04\e{-04}$ &$1.24\e{-02}$ &$   24.75
$ \\\addlinespace[.125em]
 & 
$0.25 - 0.50 $ & $7.78\e{-03}$ &$5.33\e{-04}$ &$2.58\e{-02}$ &$   48.36
$ \\\addlinespace[.125em]
 & 
$0.50 - 0.75 $ & $5.49\e{-03}$ &$3.32\e{-04}$ &$3.63\e{-02}$ &$  109.16
$ \\\addlinespace[.125em]
 & 
$0.75 - 1.00 $ & $3.84\e{-03}$ &$1.28\e{-04}$ &$4.91\e{-02}$ &$  381.27
$ \\\addlinespace[.125em]
\addlinespace[.25em]
100,000 to 500,000                                                                                   &$0.00 - 0.25 $ & $7.07\e{-03}$ &$5.73\e{-04}$ &$5.40\e{-02}$ &$   94.38
$ \\\addlinespace[.125em]
 & 
$0.25 - 0.50 $ & $1.27\e{-02}$ &$1.14\e{-03}$ &$3.14\e{-01}$ &$  275.27
$ \\\addlinespace[.125em]
 & 
$0.50 - 0.75 $ & $1.21\e{-02}$ &$9.41\e{-04}$ &$5.47\e{-01}$ &$  582.02
$ \\\addlinespace[.125em]
 & 
$0.75 - 1.00 $ & $9.30\e{-03}$ &$4.18\e{-04}$ &$7.39\e{-01}$ &$ 1766.48
$ \\\addlinespace[.125em]
\addlinespace[.25em]
500,000 to 1,000,000                                                                                 &$0.00 - 0.25 $ & $1.06\e{-02}$ &$9.22\e{-04}$ &$2.96\e{-01}$ &$  321.57
$ \\\addlinespace[.125em]
 & 
$0.25 - 0.50 $ & $1.86\e{-02}$ &$1.72\e{-03}$ &$1.07\e{+00} $ &$  619.20
$ \\\addlinespace[.125em]
 & 
$0.50 - 0.75 $ & $1.54\e{-02}$ &$1.19\e{-03}$ &$1.56\e{+00} $ &$ 1306.14
$ \\\addlinespace[.125em]
 & 
$0.75 - 1.00 $ & $1.20\e{-02}$ &$4.99\e{-04}$ &$1.98\e{+00} $ &$ 3979.82
$ \\\addlinespace[.125em]
\addlinespace[.25em]
\bottomrule
\end{tabular}

\end{center}
\caption{The average time (in seconds) required by the accelerated and unaccelerated versions
of the algorithm of this paper to construct the phase function
$\PsiP{n}(x;\gamma)$ they use to represent $\PS{n}(z;\gamma)$,
as well as the time required by the Osipov-Xiao-Rokhlin method to construct
the Legendre expansion it uses to represent $\PS{n}(z;\gamma)$.
The ``Ratio'' column reports the ratio of the time
taken by the Osipov-Xiao-Rokhlin method to the time taken by the accelerated version
of the algorithm of this paper.}
\label{experiments:table3}
\end{table}

Table~\ref{experiments:table3}  gives the average time required to construct the phase function $\PsiP{n}(x;\gamma)$
using the accelerated and unaccelerated versions of this algorithm
as well as the average time taken by the Osipov-Xiao-Rokhlin algorithm to construct
the Legendre expansion representing $\PS{n}(z;\gamma)$
for various ranges of the parameters $\gamma$  and $\sigma$.  
To generate the data for each row of the table,
first   $100$ equispaced values of the each of the parameters
$\gamma$ and $\sigma$ were sampled.  
Then, for each of the  $10,000$ pairs of sampled values of the parameter,
the phase function $\PsiP{n}(x;\gamma)$ representing $\PS{n}(z;\gamma)$ was constructed using both
the accelerated and unaccelerated algorithms of this paper and a Legendre expansion representing $\PS{n}(z;\gamma)$
was constructed with the  Osipov-Xiao-Rokhlin method.  The average time taken by each
of these procedures is reported. The column labelled `Ratio' in Table~\ref{experiments:table3} gives the
ratio of the average time taken by the Osipov-Xiao-Rokhlin method
to the average time taken by the accelerated version of our algorithm.

Each of the plots on the left side of Figure~\ref{experiments:figure3} gives the time 
required by the accelerated version of our algorithm
to construct $\PsiP{n}(x;\gamma)$ as a function of $\gamma$ for several different
values of $\sigma$.  Similarly, each plot on the right side of  Figure~\ref{experiments:figure3}
gives the time needed to construct $\PsiP{n}(x;\gamma)$ as a function of $\sigma$ for various values of $\gamma$.

Figure~\ref{experiments:figure4} contains plots comparing the time required to construct $\PsiP{n}(x;\gamma)$
using the accelerated version of the method of this paper with the time required by the Osipov-Xiao-Rokhlin algorithm
to construct the Legendre expansion representing $\PS{n}(z;\gamma)$.
Each plot on the left side of  Figure~\ref{experiments:figure4} gives these quantities as functions of
$\gamma$ for a fixed $\sigma$, while the plots on the right side of the figure
gives these quantities as functions of $\sigma$ for a fixed $\gamma$.  We choose the range
of $\gamma$ displayed in each plot on the right in order to emphasize the break-even point between
the two methods.  

From these experiments, we conclude that the time required to construct $\PsiP{n}(x;\gamma)$
grows sublogarithmically with $\gamma$ and is bounded independent of $\sigma$.  Moreover, we see that the Osipov-Xiao-Rokhlin
algorithm is faster for small values of the parameter, but for  $\gamma \sigma > 250$  or so, the accelerated
version of the algorithm of this paper becomes more efficient, and it is much more efficient
at large values of $\gamma$.

\end{subsection}

\begin{subsection}{The time required to evaluate $\PS{n}(z;\gamma)$}

In this final set of experiments, we measured the time required to evaluate $\PS{n}(z;\gamma)$
using the phase function $\PsiP{n}(x;\gamma)$ and compared it to the 
time needed to evaluate the Legendre expansion used by the Osipov-Xiao-Rokhlin to represent 
$\PS{n}(z;\gamma)$.    Table~\ref{experiments:table4} and Figures~\ref{experiments:figure5}  
give the results.

Table~\ref{experiments:table4}  gives the average time required to evaluate $\PS{n}(z;\gamma)$
using the the phase function $\PsiP{n}(x;\gamma)$
as well as the average time required to evaluate the Legendre expansion used by the Osipov-Xiao-Rokhlin algorithm to 
represent $\PS{n}(z;\gamma)$  for various ranges of the parameters $\gamma$  and $\sigma$.  
To generate the data for each row of the table,
we first sampled  $20$ equispaced values of the each of the parameters
$\gamma$ and $\sigma$.
Then, for each of the  $400$ pairs of sampled values of the parameter,
the phase function $\PsiP{n}(x;\gamma)$ representing $\PS{n}(z;\gamma)$ was constructed using
the algorithm of this paper and a Legendre expansion representing $\PS{n}(z;\gamma)$
was constructed with the  Osipov-Xiao-Rokhlin method.    Next, for each pair of the parameters
$\gamma$ and $n$, $\PS{n}(z;\gamma)$ was evaluated at $100$ equispaced points on the interval $(0,1)$.
The average time taken to perform these $40,000$ evaluations is reported.

Each of the plots on the left side of Figure~\ref{experiments:figure5} gives the time 
required to evaluate $\PsiP{n}(x;\gamma)$ as a function of $\gamma$ for several different
values of $\sigma$.  Similarly, each plot on the right side of  Figure~\ref{experiments:figure3}
gives the time needed to evaluate $\PS{n}(z;\gamma)$ as a function of $\sigma$ for various values of $\gamma$.

We omit the analog of Figure~\ref{experiments:figure4} comparing the time required to evaluate
the phase function $\PsiP{n}(x;\gamma)$ used by the algorithm of this paper to represent
$\PS{n}(z;\gamma)$  and the time required to evaluate it using the 
Legendre expansion used by Xiao-Osipov-Rokhlin method to represent $\PS{n}(z;\gamma)$.
This is because, even for small values of $\gamma$ and $\sigma$, the 
phase function method is considerably faster.  

We conclude from these experiments that the time required to evaluate $\PS{n}(z;\gamma)$
using the phase function $\PsiP{n}(x;\gamma)$ is largely independent of the parameters
$\gamma$ and $\sigma$.

\begin{table}[!h]
\begin{center}
\begin{tabular}{cccccr}
\toprule
Range of $\gamma$ & Range of $\sigma$ &   Average time & Average time        & Ratio\\
              &                   & Phase  & Rokhlin, et. al.\\
\midrule
100 to 500                                                                                           &$0.00 - 0.25 $ & $1.22\e{-07}$ &$1.66\e{-06}$ &$     13.67
$ \\\addlinespace[.125em]
 & 
$0.25 - 0.50 $ & $1.14\e{-07}$ &$2.17\e{-06}$ &$     18.99
$ \\\addlinespace[.125em]
 & 
$0.50 - 0.75 $ & $1.36\e{-07}$ &$2.44\e{-06}$ &$     17.88
$ \\\addlinespace[.125em]
 & 
$0.75 - 1.00 $ & $1.16\e{-07}$ &$2.70\e{-06}$ &$     23.30
$ \\\addlinespace[.125em]
\addlinespace[.25em]
500 to 1,000                                                                                         &$0.00 - 0.25 $ & $1.36\e{-07}$ &$3.28\e{-06}$ &$     24.08
$ \\\addlinespace[.125em]
 & 
$0.25 - 0.50 $ & $1.18\e{-07}$ &$4.59\e{-06}$ &$     38.74
$ \\\addlinespace[.125em]
 & 
$0.50 - 0.75 $ & $1.17\e{-07}$ &$5.29\e{-06}$ &$     45.25
$ \\\addlinespace[.125em]
 & 
$0.75 - 1.00 $ & $1.18\e{-07}$ &$5.89\e{-06}$ &$     49.92
$ \\\addlinespace[.125em]
\addlinespace[.25em]
1,000 to 5,000                                                                                       &$0.00 - 0.25 $ & $1.55\e{-07}$ &$9.89\e{-06}$ &$     63.71
$ \\\addlinespace[.125em]
 & 
$0.25 - 0.50 $ & $1.26\e{-07}$ &$1.57\e{-05}$ &$    124.06
$ \\\addlinespace[.125em]
 & 
$0.50 - 0.75 $ & $1.17\e{-07}$ &$1.86\e{-05}$ &$    158.33
$ \\\addlinespace[.125em]
 & 
$0.75 - 1.00 $ & $1.17\e{-07}$ &$2.10\e{-05}$ &$    178.54
$ \\\addlinespace[.125em]
\addlinespace[.25em]
5,000 to 10,000                                                                                      &$0.00 - 0.25 $ & $1.62\e{-07}$ &$2.23\e{-05}$ &$    137.12
$ \\\addlinespace[.125em]
 & 
$0.25 - 0.50 $ & $1.30\e{-07}$ &$3.73\e{-05}$ &$    287.03
$ \\\addlinespace[.125em]
 & 
$0.50 - 0.75 $ & $1.18\e{-07}$ &$4.49\e{-05}$ &$    378.10
$ \\\addlinespace[.125em]
 & 
$0.75 - 1.00 $ & $1.18\e{-07}$ &$5.09\e{-05}$ &$    429.64
$ \\\addlinespace[.125em]
\addlinespace[.25em]
10,000 to 50,000                                                                                     &$0.00 - 0.25 $ & $1.66\e{-07}$ &$8.18\e{-05}$ &$    491.98
$ \\\addlinespace[.125em]
 & 
$0.25 - 0.50 $ & $1.32\e{-07}$ &$1.43\e{-04}$ &$   1088.13
$ \\\addlinespace[.125em]
 & 
$0.50 - 0.75 $ & $1.19\e{-07}$ &$1.74\e{-04}$ &$   1465.68
$ \\\addlinespace[.125em]
 & 
$0.75 - 1.00 $ & $1.17\e{-07}$ &$1.98\e{-04}$ &$   1683.48
$ \\\addlinespace[.125em]
\addlinespace[.25em]
50,000 to 100,000                                                                                    &$0.00 - 0.25 $ & $1.66\e{-07}$ &$1.99\e{-04}$ &$   1196.92
$ \\\addlinespace[.125em]
 & 
$0.25 - 0.50 $ & $1.35\e{-07}$ &$3.58\e{-04}$ &$   2645.23
$ \\\addlinespace[.125em]
 & 
$0.50 - 0.75 $ & $1.21\e{-07}$ &$4.35\e{-04}$ &$   3597.88
$ \\\addlinespace[.125em]
 & 
$0.75 - 1.00 $ & $1.20\e{-07}$ &$4.96\e{-04}$ &$   4129.21
$ \\\addlinespace[.125em]
\addlinespace[.25em]
100,000 to 500,000                                                                                   &$0.00 - 0.25 $ & $1.61\e{-07}$ &$7.87\e{-04}$ &$   4864.34
$ \\\addlinespace[.125em]
 & 
$0.25 - 0.50 $ & $1.39\e{-07}$ &$1.61\e{-03}$ &$  11542.75
$ \\\addlinespace[.125em]
 & 
$0.50 - 0.75 $ & $1.33\e{-07}$ &$2.75\e{-03}$ &$  20682.52
$ \\\addlinespace[.125em]
 & 
$0.75 - 1.00 $ & $1.35\e{-07}$ &$3.27\e{-03}$ &$  24266.13
$ \\\addlinespace[.125em]
\addlinespace[.25em]
500,000 to 1,000,000                                                                                 &$0.00 - 0.25 $ & $1.70\e{-07}$ &$2.49\e{-03}$ &$  14689.04
$ \\\addlinespace[.125em]
 & 
$0.25 - 0.50 $ & $1.78\e{-07}$ &$1.28\e{-02}$ &$  72054.88
$ \\\addlinespace[.125em]
 & 
$0.50 - 0.75 $ & $1.64\e{-07}$ &$1.73\e{-02}$ &$ 105582.86
$ \\\addlinespace[.125em]
 & 
$0.75 - 1.00 $ & $1.58\e{-07}$ &$1.80\e{-02}$ &$ 114213.05
$ \\\addlinespace[.125em]
\addlinespace[.25em]
\bottomrule
\end{tabular}

\end{center}
\caption{The average time (in seconds) required by the algorithm of this paper and by the 
Osipov-Xiao-Rokhlin method to evaluate $\PS{n}(z;\gamma)$.
The ``Ratio'' column reports the ratio of the time taken by the two methods.}
\label{experiments:table4}
\end{table}

\end{subsection}

\label{section:experiments}
\end{section}

\begin{section}{Conclusions}

It is well known that many second order differential equations
have solutions whose logarithms are easier to represent than the solutions themselves.
Historically, this observation has mainly been used to construct asymptotic
expansions of the solutions of such equations.
Here we suggest a more direct approach: the numerical computation of the
logarithms.  We discuss the results of experiments which indicate such an
approach leads to a method for the numerical evaluation of the prolate
spheroidal wave functions of order zero whose running time grows
much more slowly with bandlimit and characteristic exponent than standard algorithms.
Moreveor, the algorithm presented here can be viewed
as a template suitable for application to many other
families of special functions satisfying second order differential equations,
such as the Jacobi polynomials and Hermite functions.

We rely on experimental evidence for our claims.  It would be of great
interest to prove a bound on the cost of representing the phase function
$\PsiP{n}(x;\gamma)$ or a related function using piecewise polynomial expansions.
Proving the conjectures made in this article regarding the associated modulus function
might provide a good starting point for such an endeavor.
Formulas which imply the monotonicity properties
of many related second order differential equations
(such as those established in  \cite{durand75} and \cite{durand78})
have been proved using Koornwinder's addition formula
for Jacobi functions \cite{koornwinder}.  An analagous result
for the spheroidal wave functions would most likely lead to proofs
of the conjectures made in Section~\ref{section:monotonicity}.

We believe it is  possible to eliminate the dependence
of the running time of the algorithm described here on the bandlimit $\gamma$.
The cost of representing the logarithms of the WKB solution in regions bounded
away from the turning points of the reduced spheroidal wave equation
appears to be independent of $\gamma$.  It is only near the turning point
where the complexity increases  with $\gamma$.  A more sophisticated
approach, which used an alternate representation of $\PS{n}(z;\gamma)$
near turning points, would most likely yield an $\mathcal{O}\left(1\right)$
algorithm for evaluation the prolate spheroidal wave functions of order zero.

Finally, we note that the methodology discussed here and in \cite{prolates2} can be
applied to the spheroidal wave functions of nonzero orders.  
We believe the conjectures of Section~\ref{section:monotonicity}
extend essentially without modification to that equation.
The only additional difficulty is that the expansions used in 
\cite{prolates2} to evaluate the Sturm-Liouville eigenvalues
of the spheroidal wave functions would depend on three variables
instead of two, and be correspondingly larger and cost more to evaluate.






\label{section:conclusion}
\end{section}

\begin{section}{Acknowledgements}
The author was supported in part by NSERC Discovery grant  RGPIN-2021-02613, and
by NSF grants DMS-1818820 and  DMS-2012487.
\end{section}

\bibliographystyle{acm}
\bibliography{prolates.bib}

\vfill\eject


\begin{figure}[!h]
\hfil
\includegraphics[width=.49\textwidth]{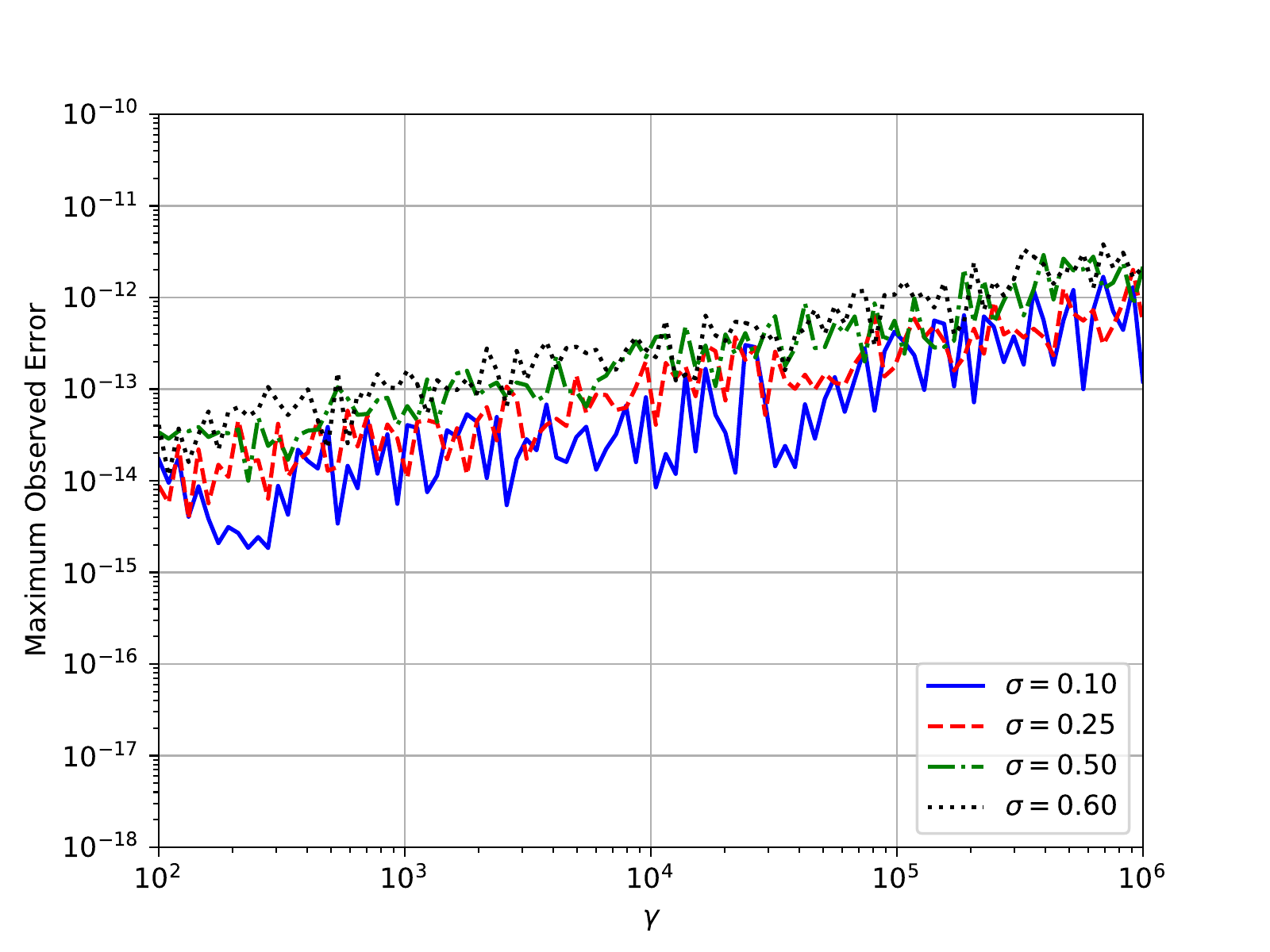}
\hfil
\includegraphics[width=.49\textwidth]{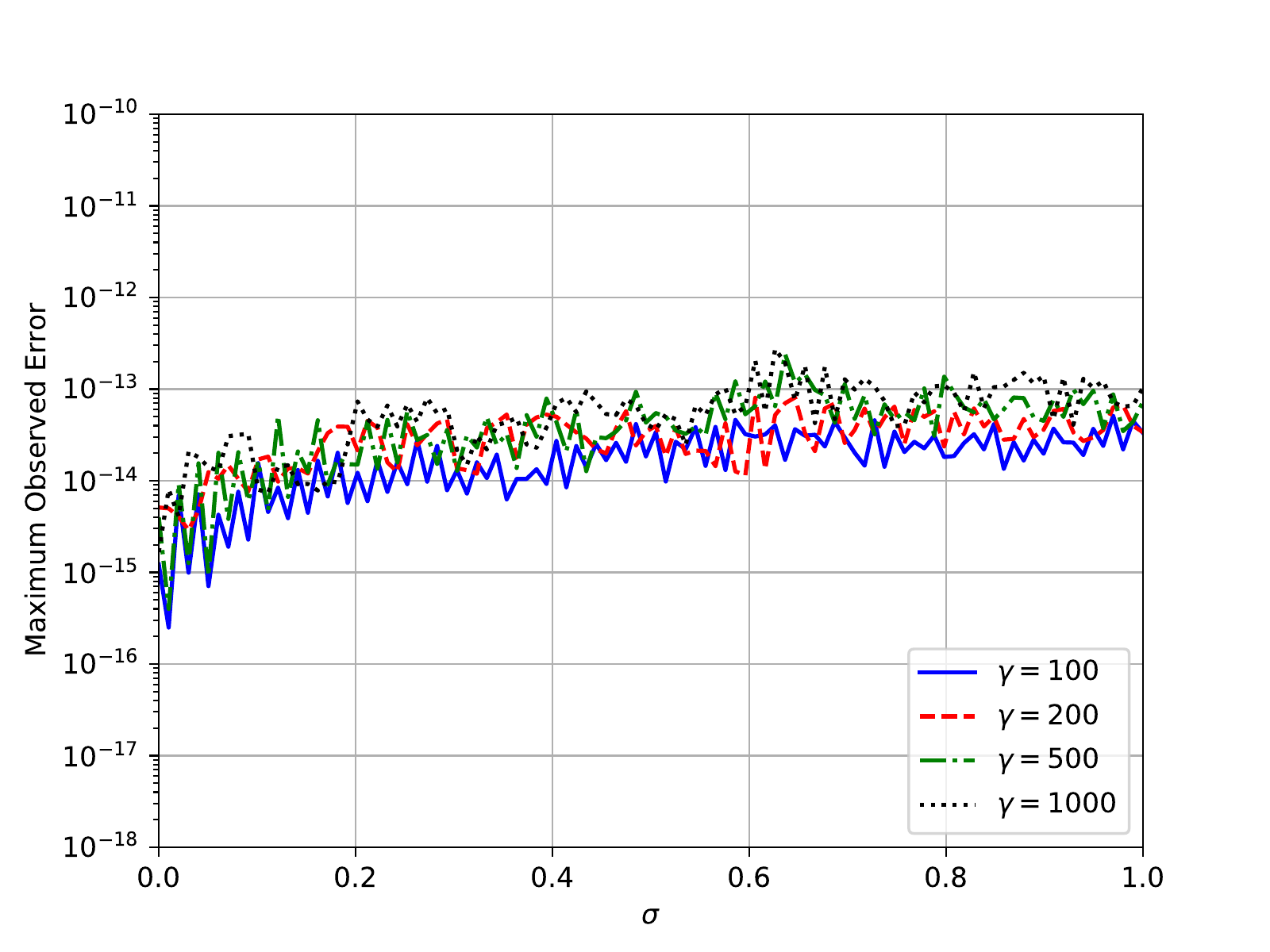}
\hfil

\vskip 3em

\hfil
\includegraphics[width=.49\textwidth]{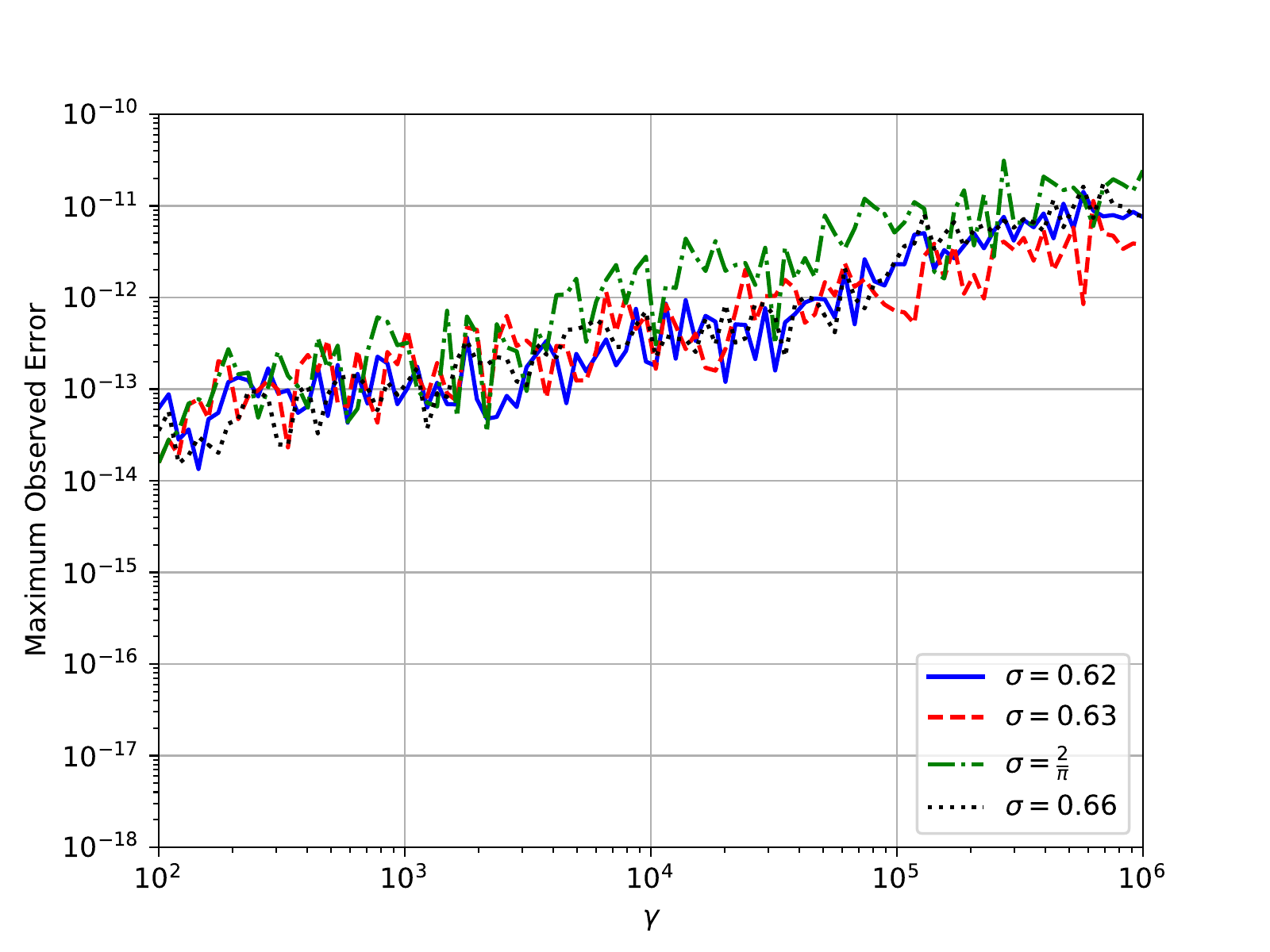}
\hfil
\includegraphics[width=.49\textwidth]{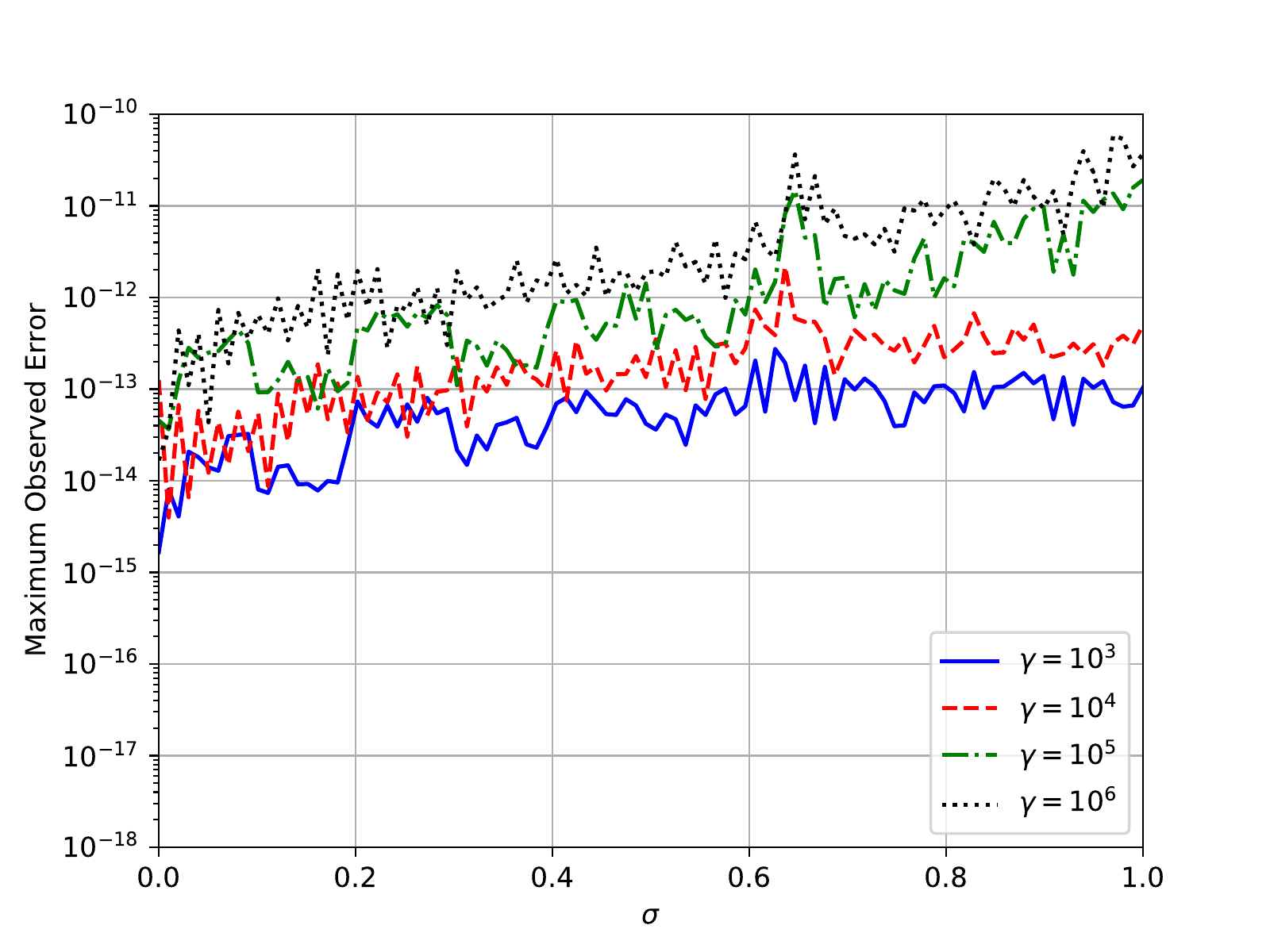}
\hfil

\vskip 3em

\hfil
\includegraphics[width=.49\textwidth]{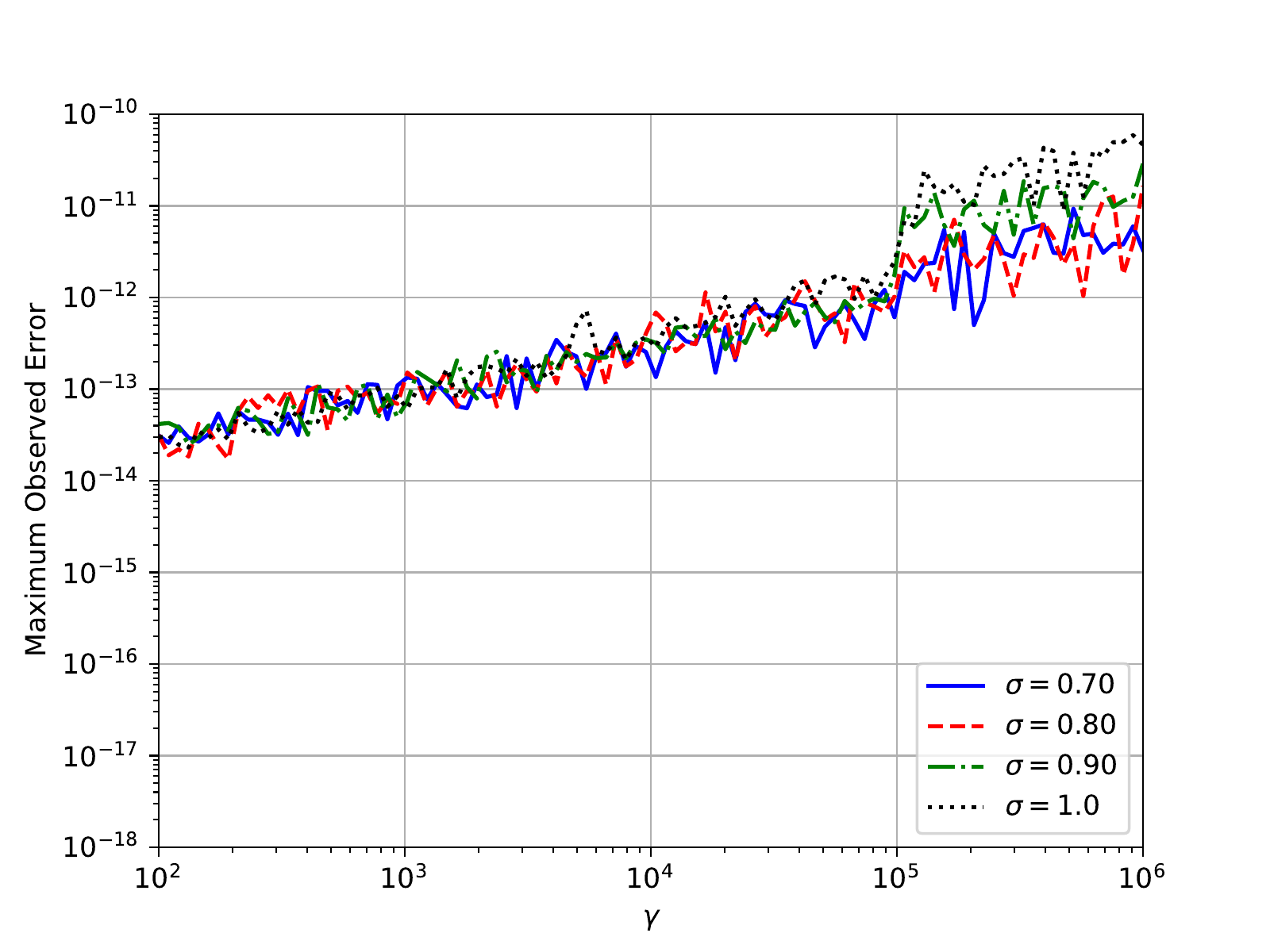}
\hfil
\includegraphics[width=.49\textwidth]{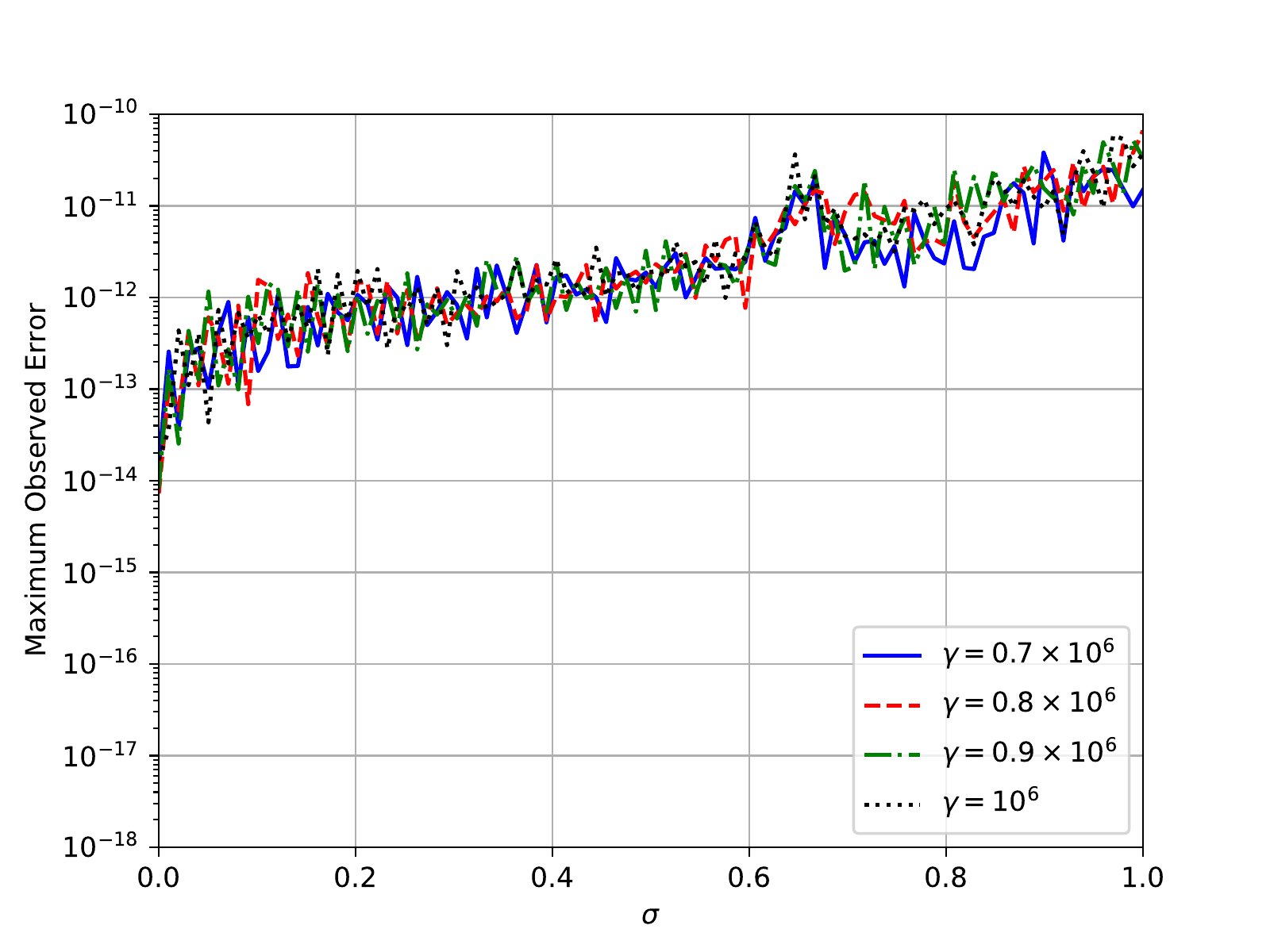}
\hfil

\caption{The accuracy with which $\PS{n}(\gamma)$ is evaluated.
Each of the graphs on the left gives the maximum observed error
as a function of $\gamma$ for several values of $\sigma$, while the plots on the right 
give the maximum observed error  as a function of $\sigma$ for several values of $\gamma$.  
A logarithmic scale is used for the x-axis in each of the plots on the left.}
\label{experiments:figure2}
\end{figure}

\begin{figure}[!t]
\hfil
\includegraphics[width=.49\textwidth]{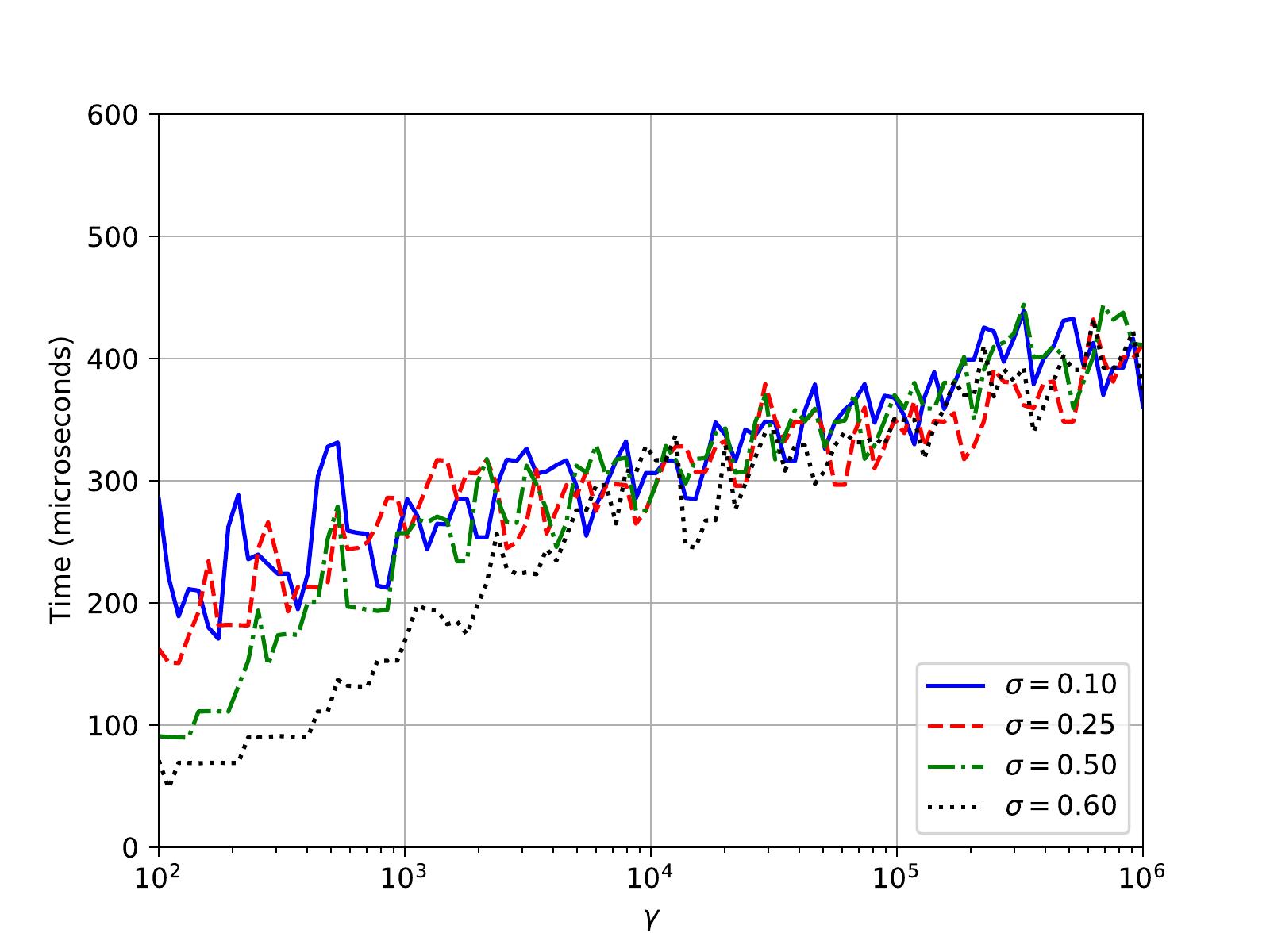}
\hfil
\includegraphics[width=.49\textwidth]{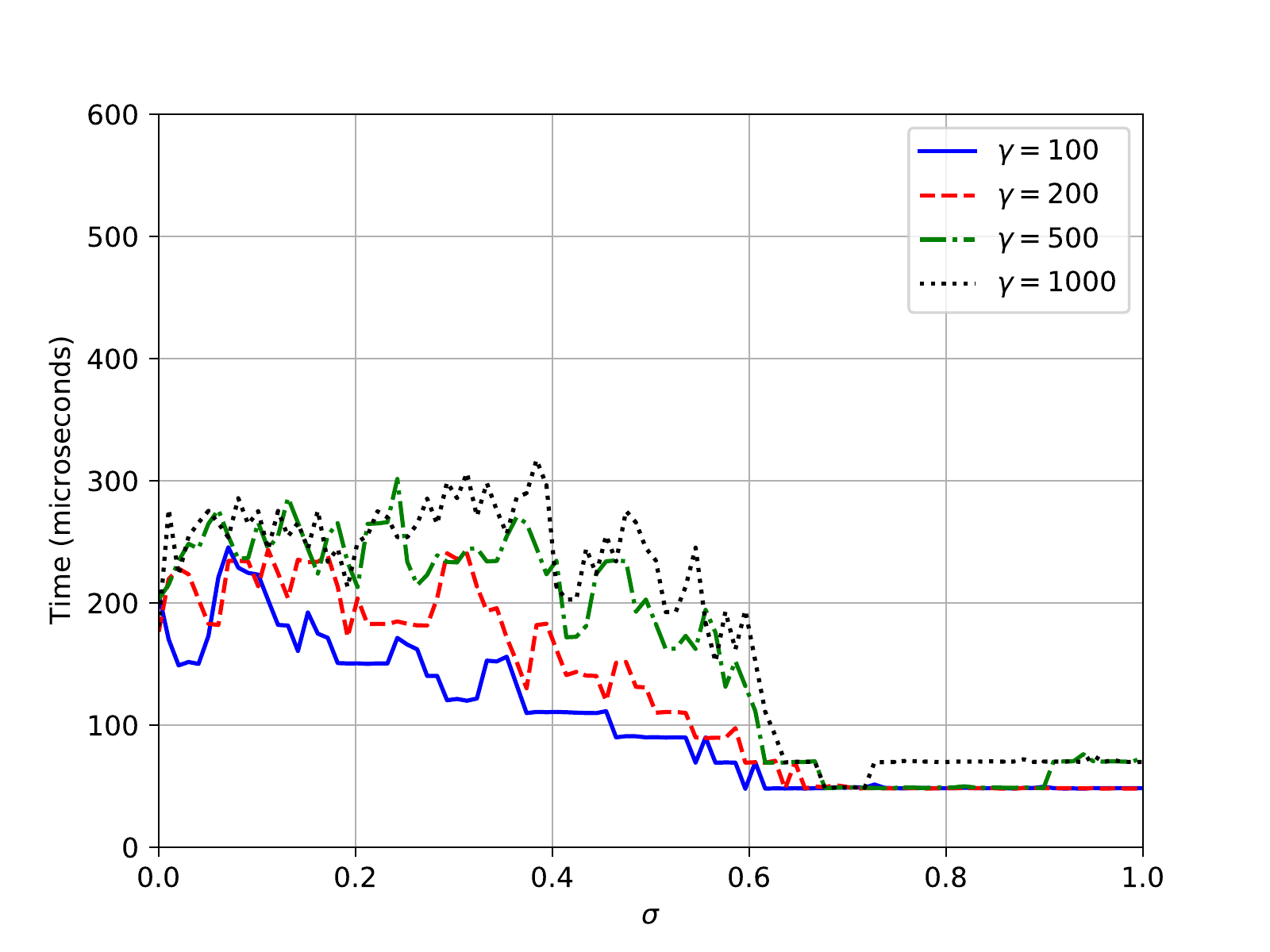}
\hfil

\vskip 3em

\hfil
\includegraphics[width=.49\textwidth]{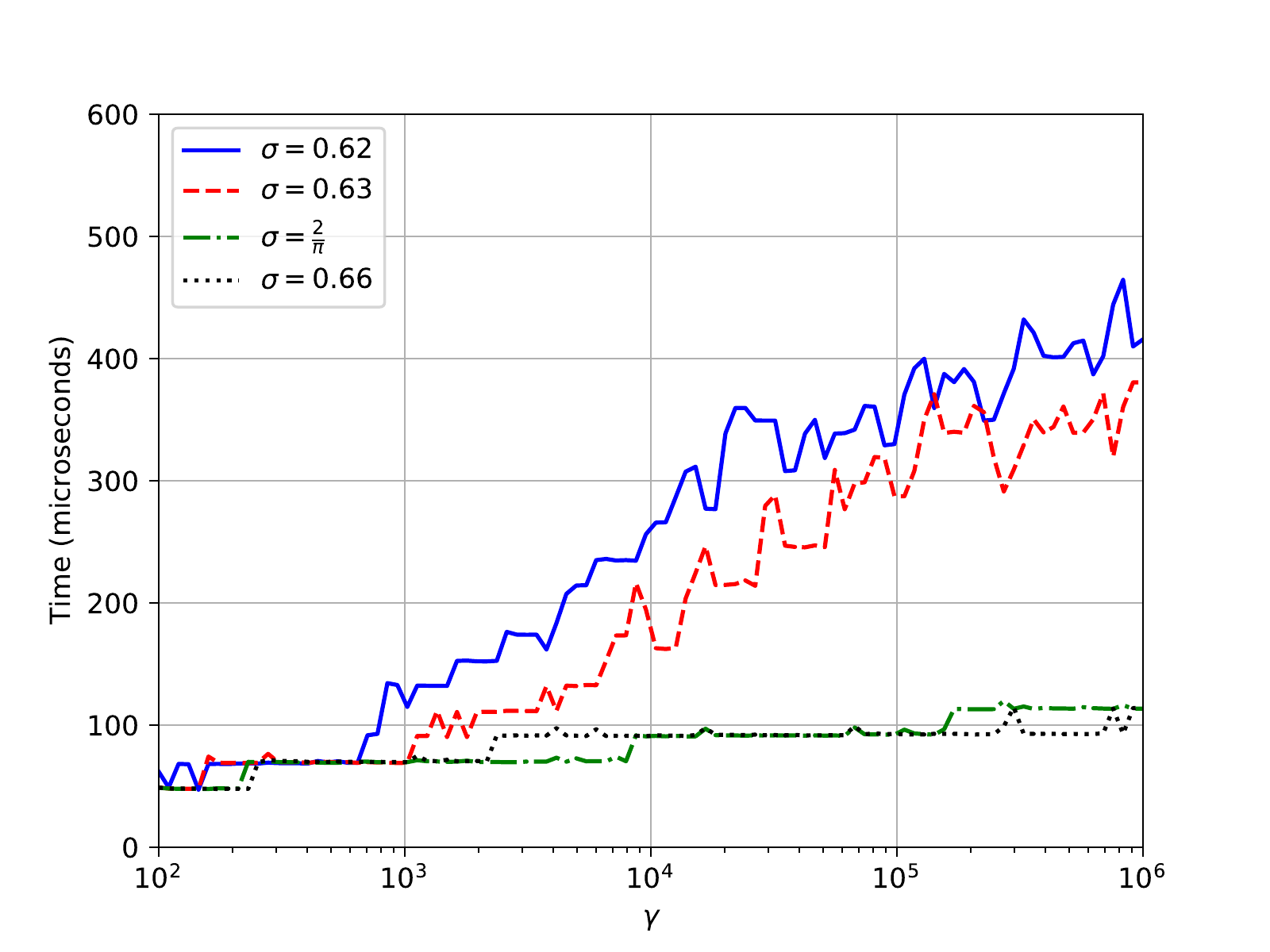}
\hfil
\includegraphics[width=.49\textwidth]{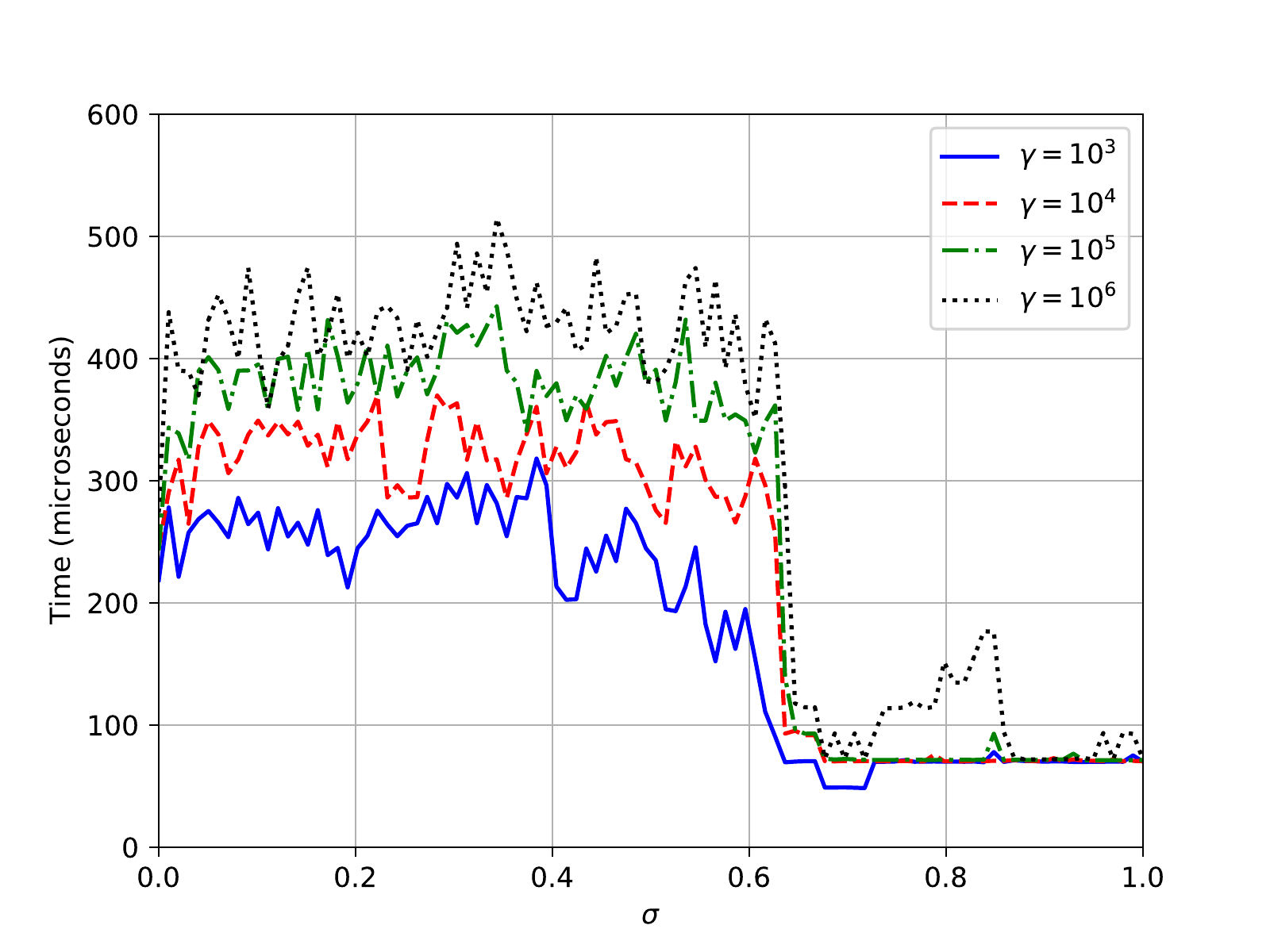}
\hfil

\vskip 3em

\hfil
\includegraphics[width=.49\textwidth]{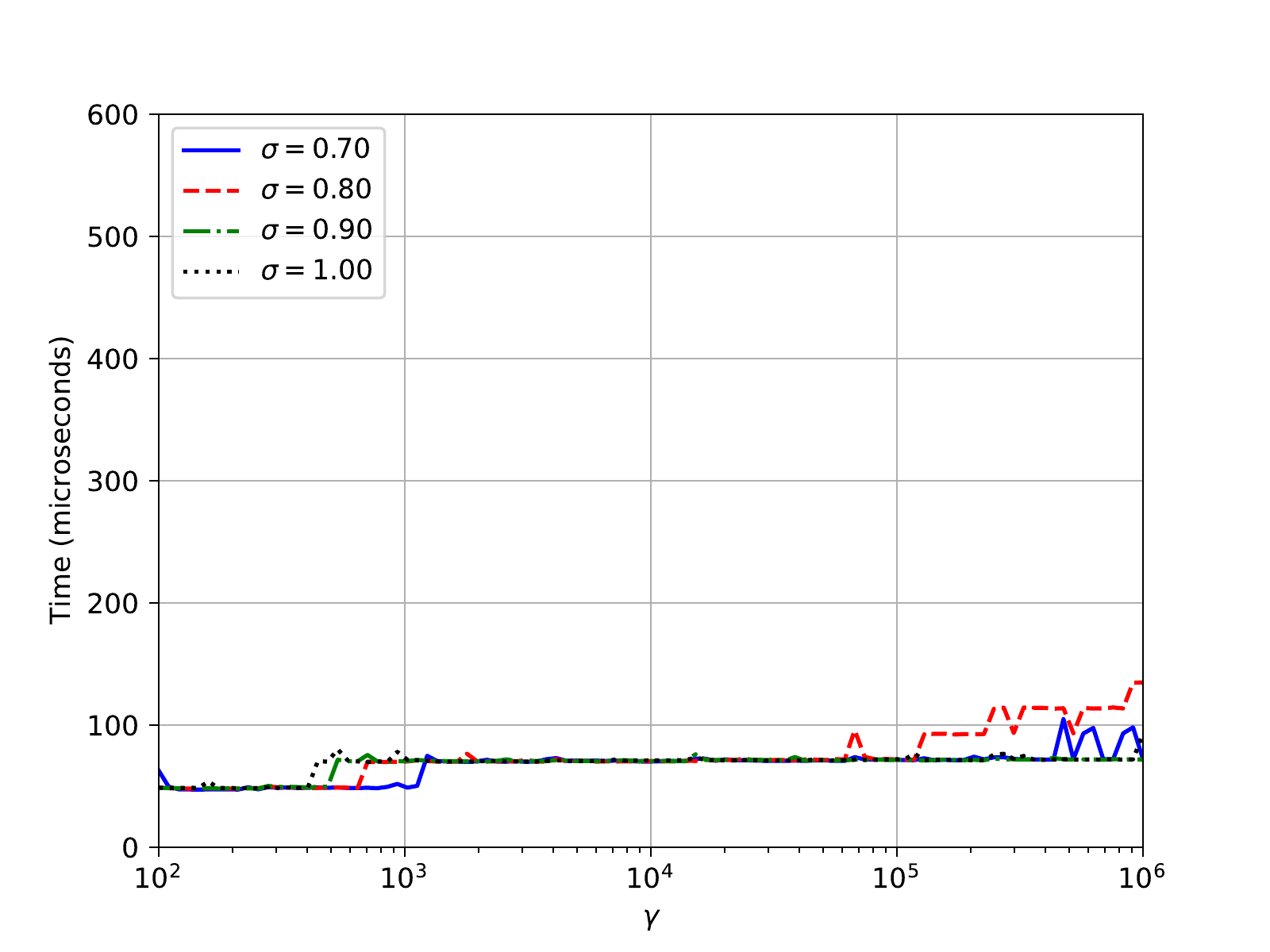}
\hfil
\includegraphics[width=.49\textwidth]{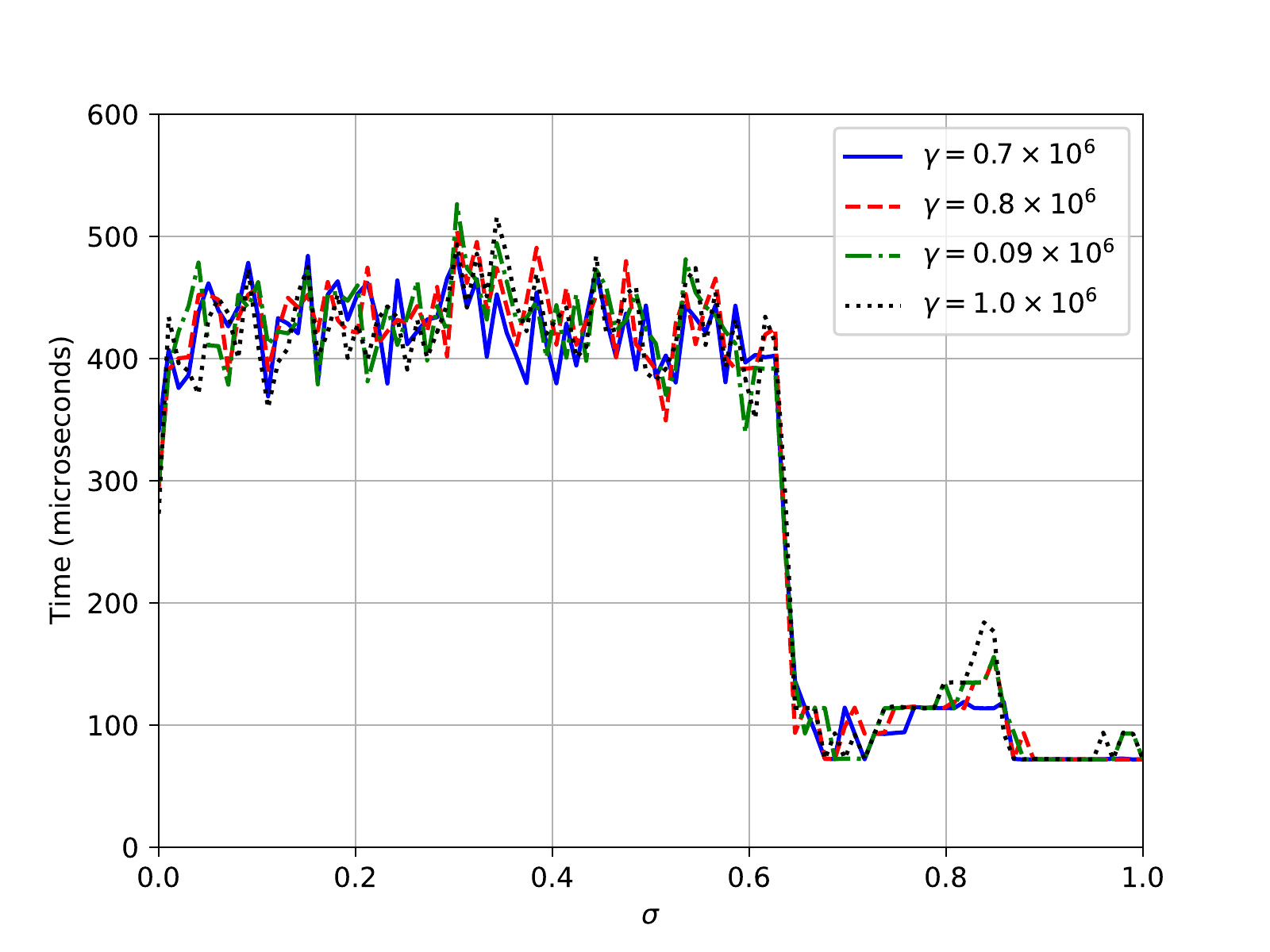}
\hfil

\caption{The time (in microseconds) required by the accelerated version of the algorithm of this paper 
to construct $\PsiP{n}(\gamma;z)$.  Each of the graphs on the left gives the time required
as a function of $\gamma$ for several values of $\sigma$, while the plots on the right 
give the required time  as a function of $\sigma$ for several values of $\gamma$. 
A logarithmic scale is used for the x-axis in each of the plots on the left.}
\label{experiments:figure3}
\end{figure}

\begin{figure}[!h]
\hfil
\includegraphics[width=.49\textwidth]{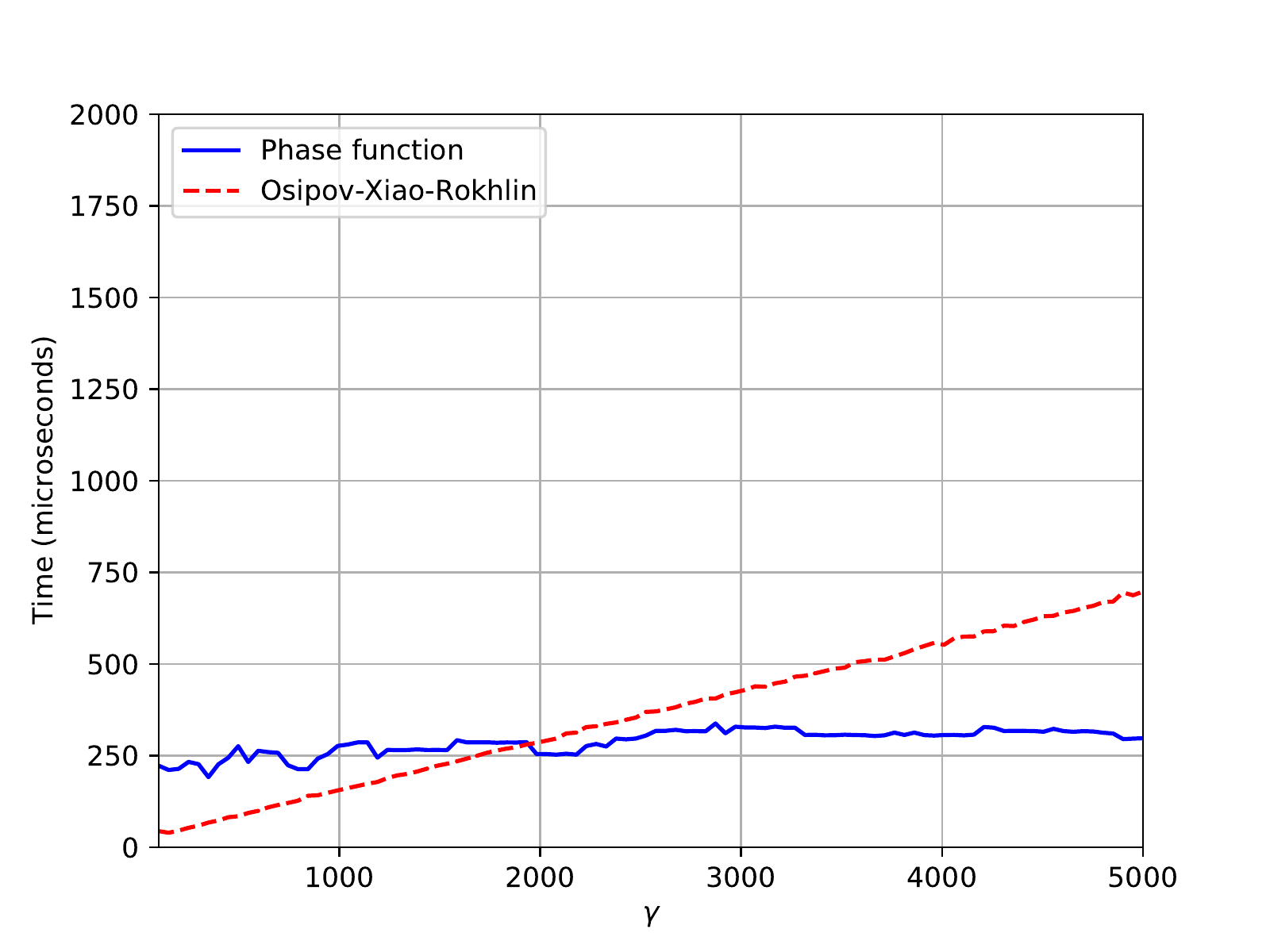}
\hfil
\includegraphics[width=.49\textwidth]{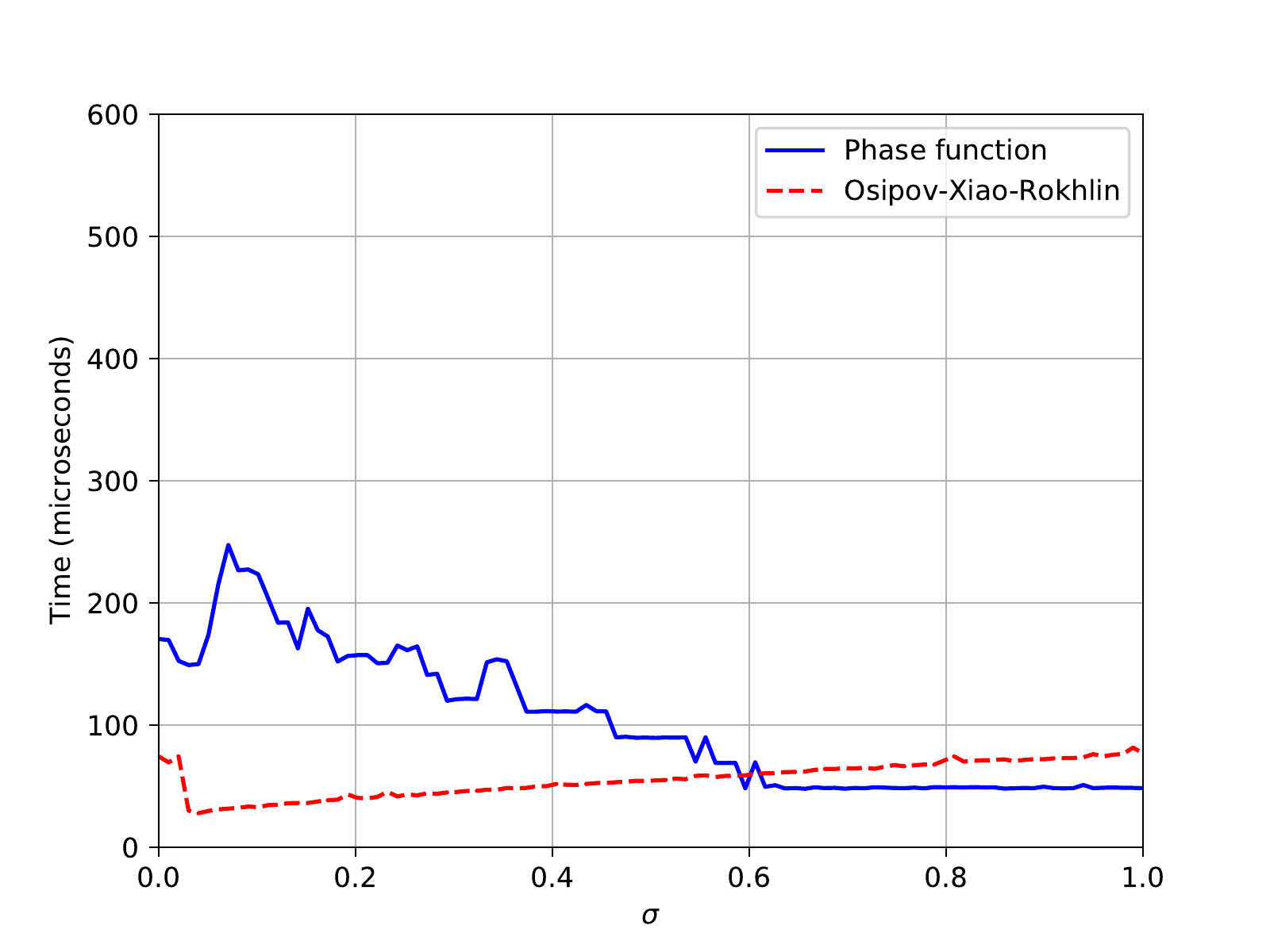}
\hfil

\vskip 3em

\hfil
\includegraphics[width=.49\textwidth]{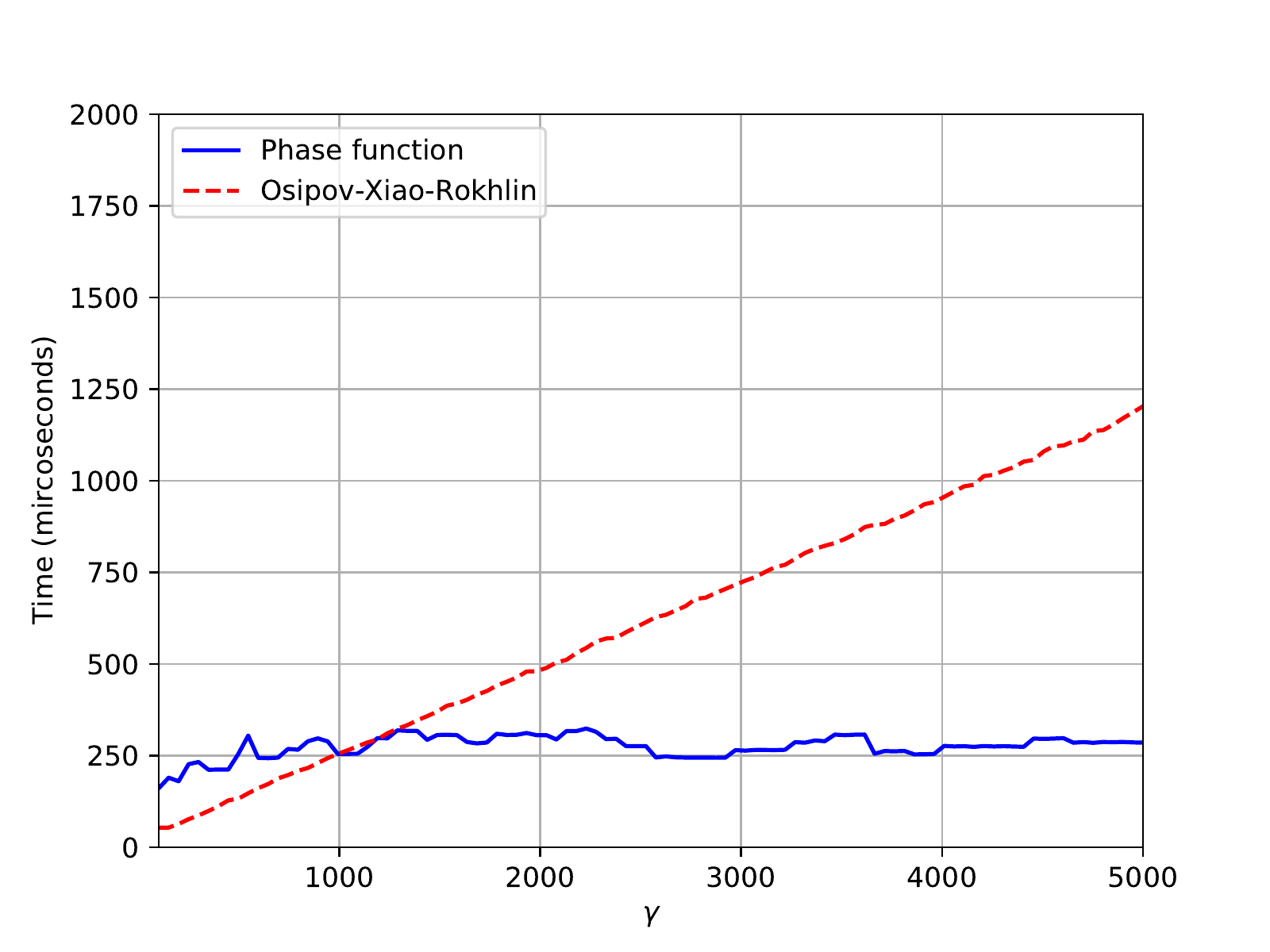}
\hfil
\includegraphics[width=.49\textwidth]{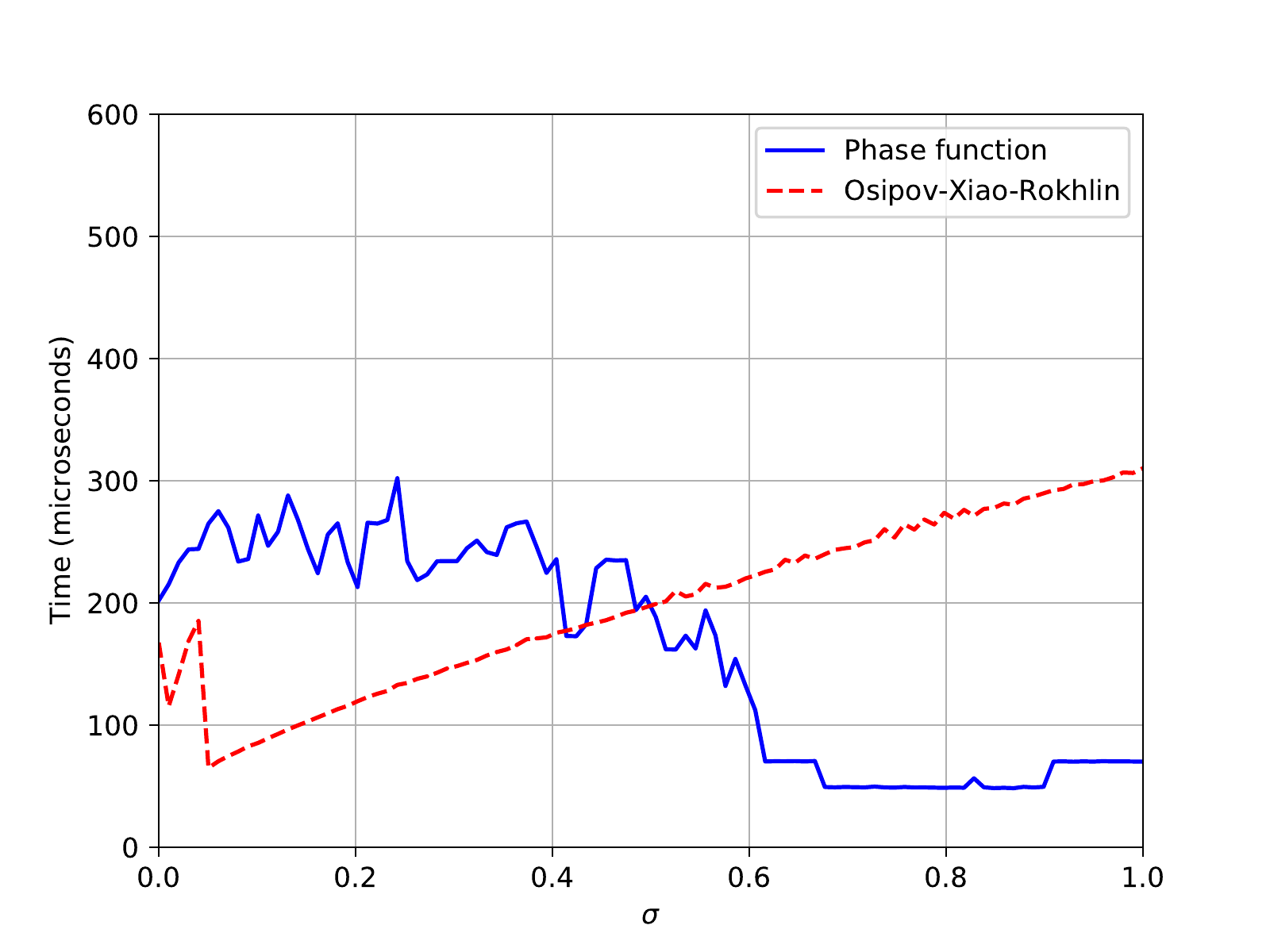}
\hfil

\vskip 3em

\hfil
\includegraphics[width=.49\textwidth]{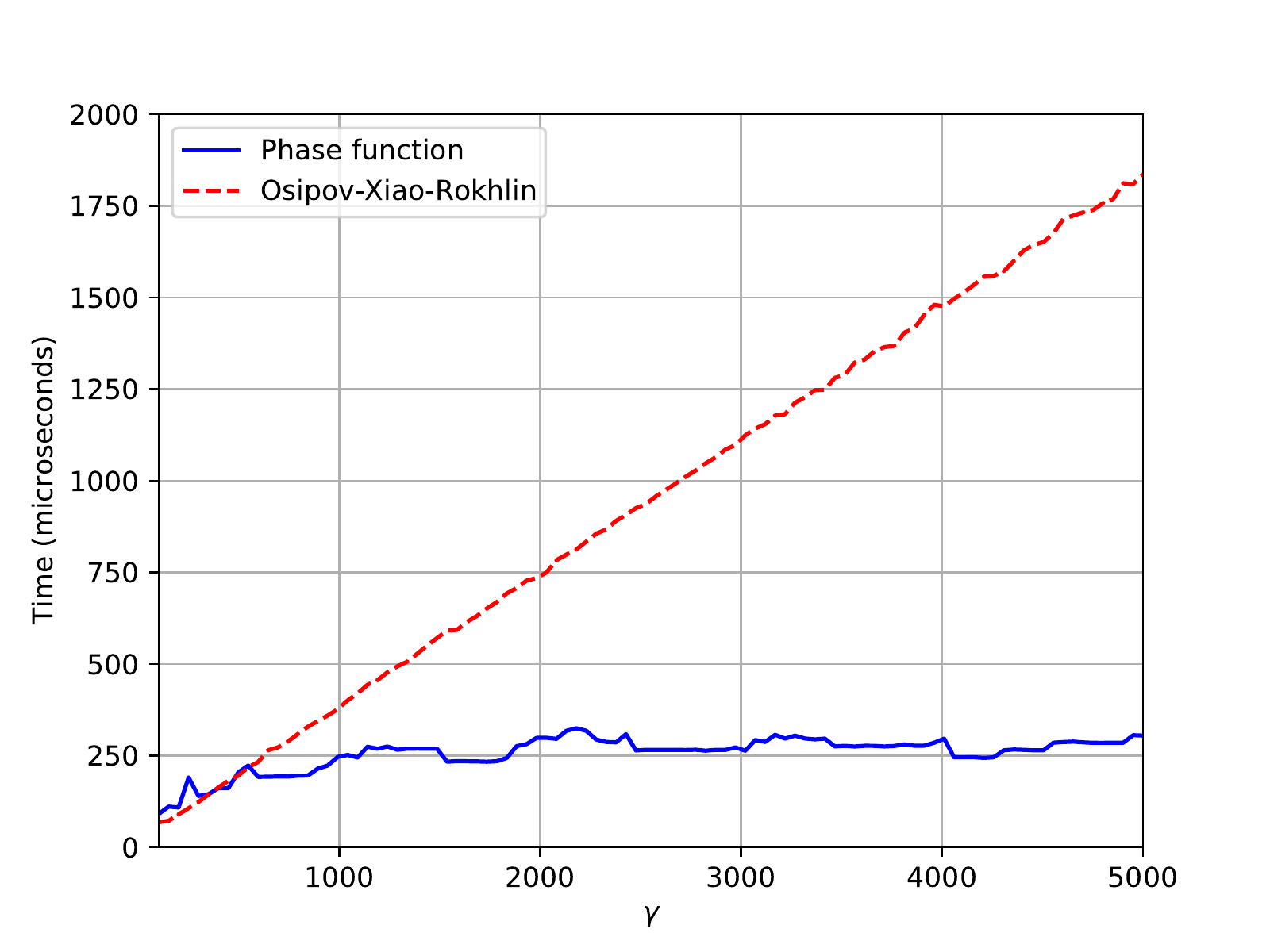}
\hfil
\includegraphics[width=.49\textwidth]{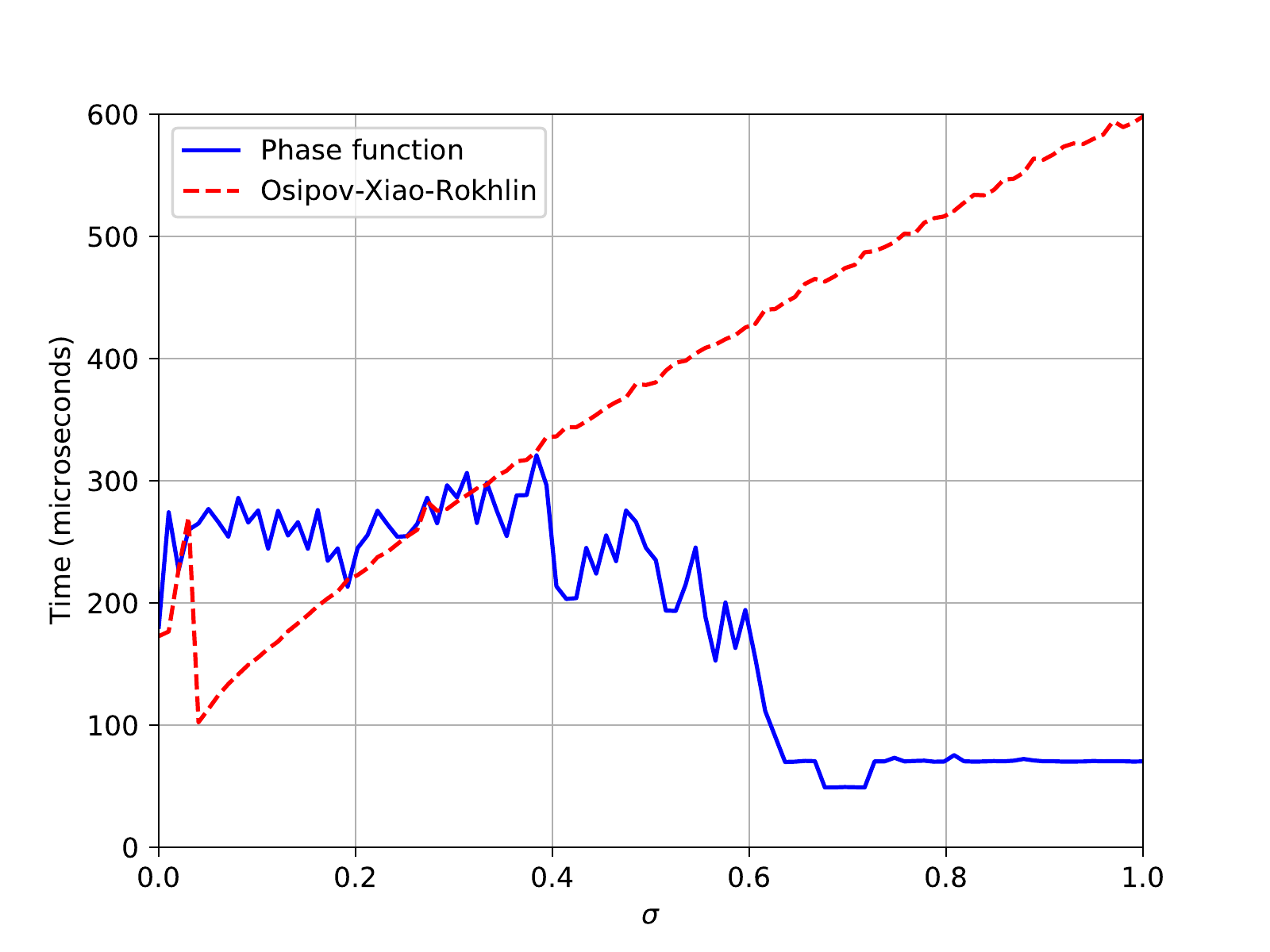}
\hfil

\caption{A comparison of the time required to construct a representation of $\PS{n}(z;\gamma)$ using
the accelerated version of the algorithm of this paper and the Xiao-Osipov-Rokhlin method.  
Each plot on the left  gives these quanitites as a function of gamma for a fixed value of $\sigma$.   
From top to bottom, $\sigma = 0.10$, $\sigma=0.25$ and $\sigma=0.50$.
Each plot on the right  gives these quanitites as a function of sigma for fixed values of $\gamma$.   
From top to bottom, $\gamma = 100$, $\gamma=500$ and $\gamma=1,000$.
}

\label{experiments:figure4}
\end{figure}

\begin{figure}[!t]
\hfil
\includegraphics[width=.49\textwidth]{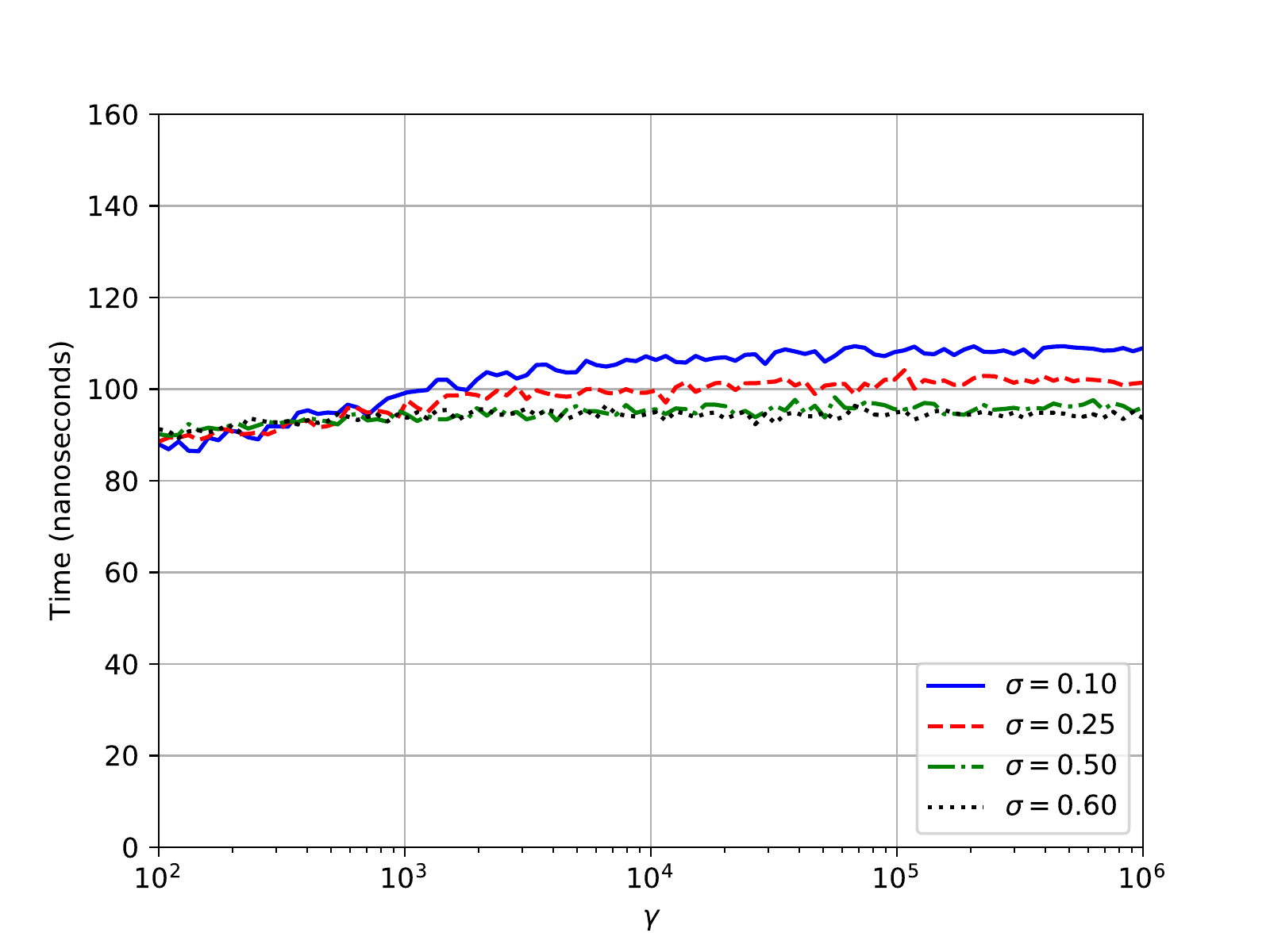}
\hfil
\includegraphics[width=.49\textwidth]{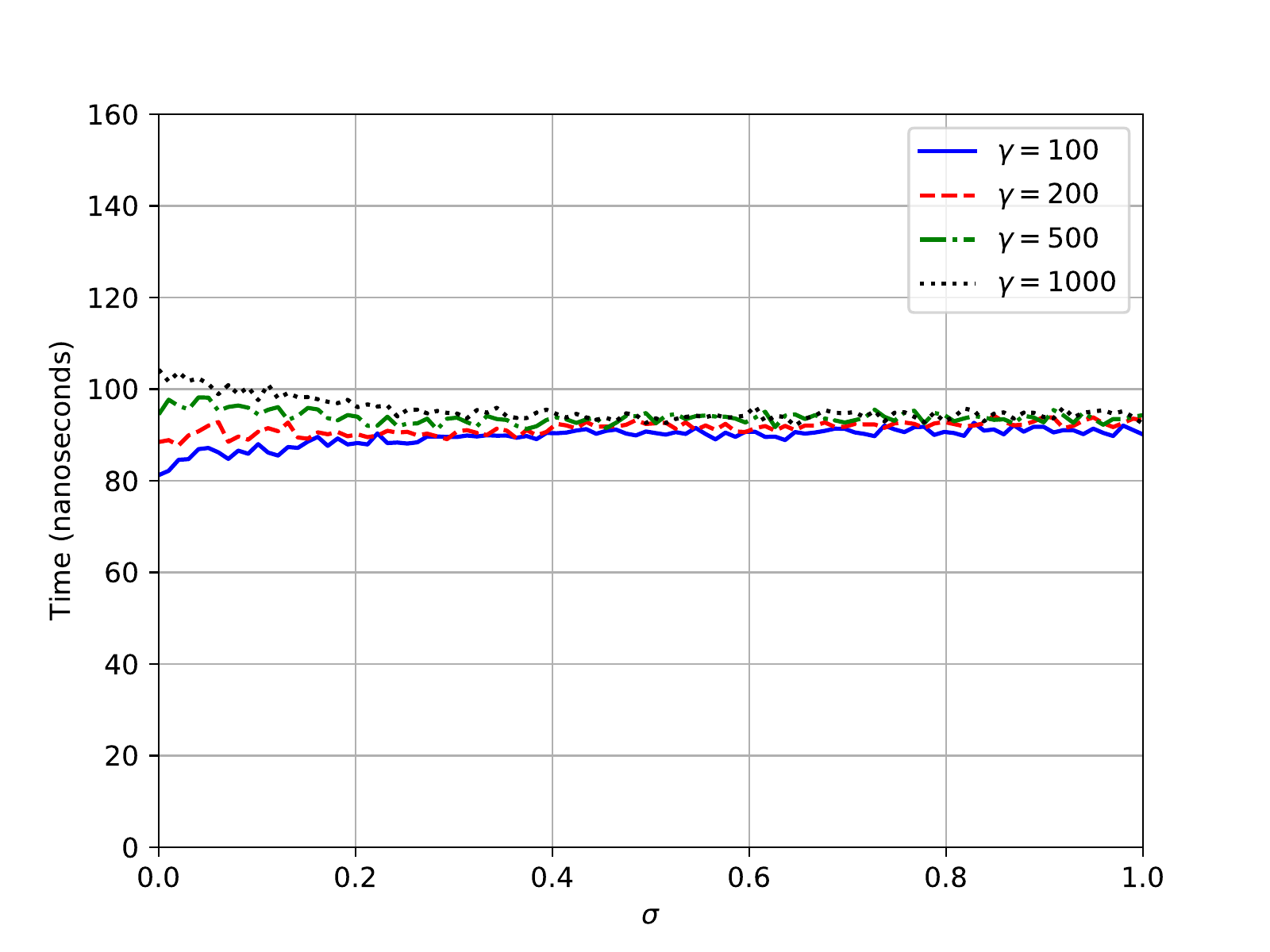}
\hfil

\vskip 3em

\hfil
\includegraphics[width=.49\textwidth]{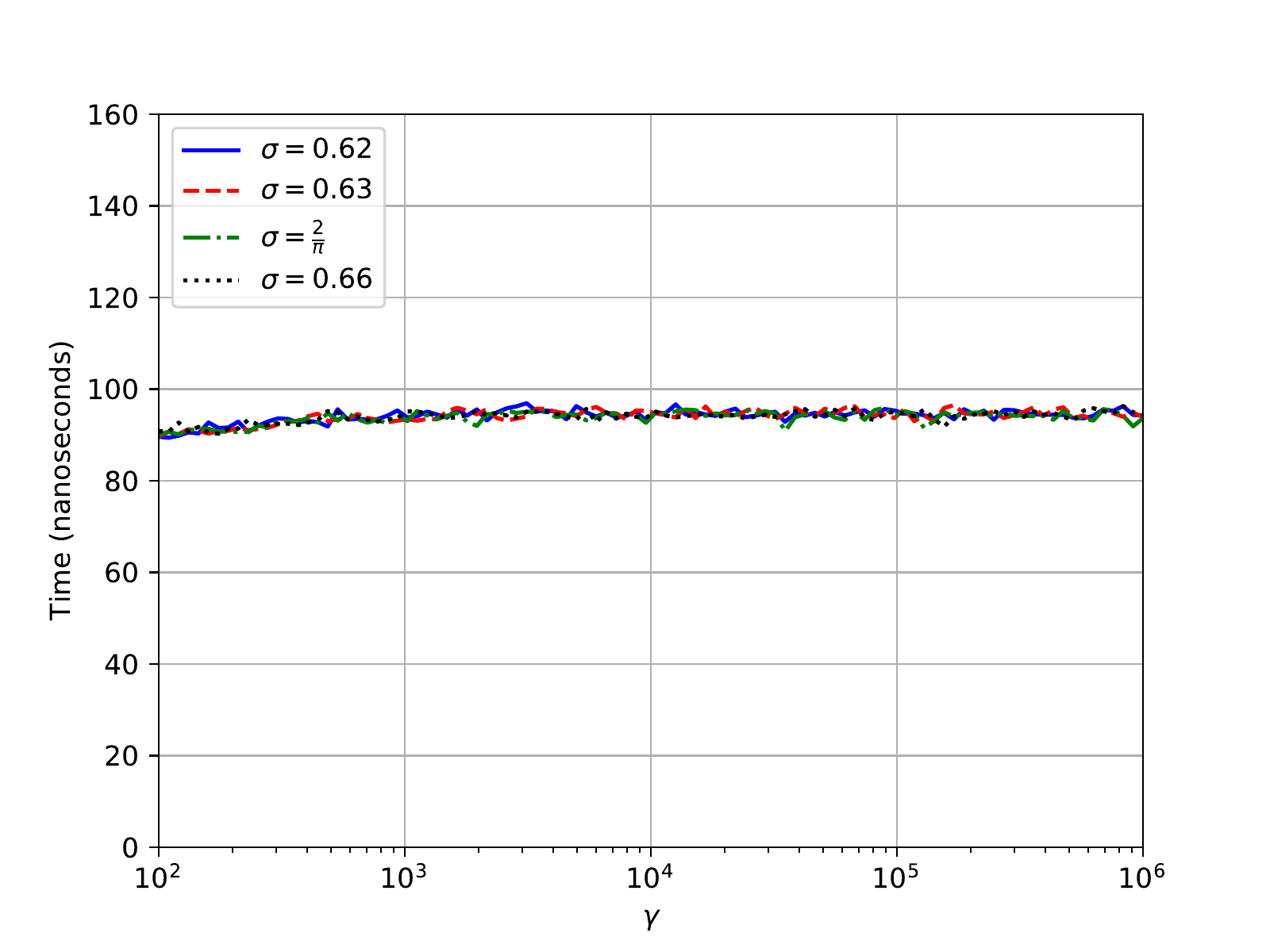}
\hfil
\includegraphics[width=.49\textwidth]{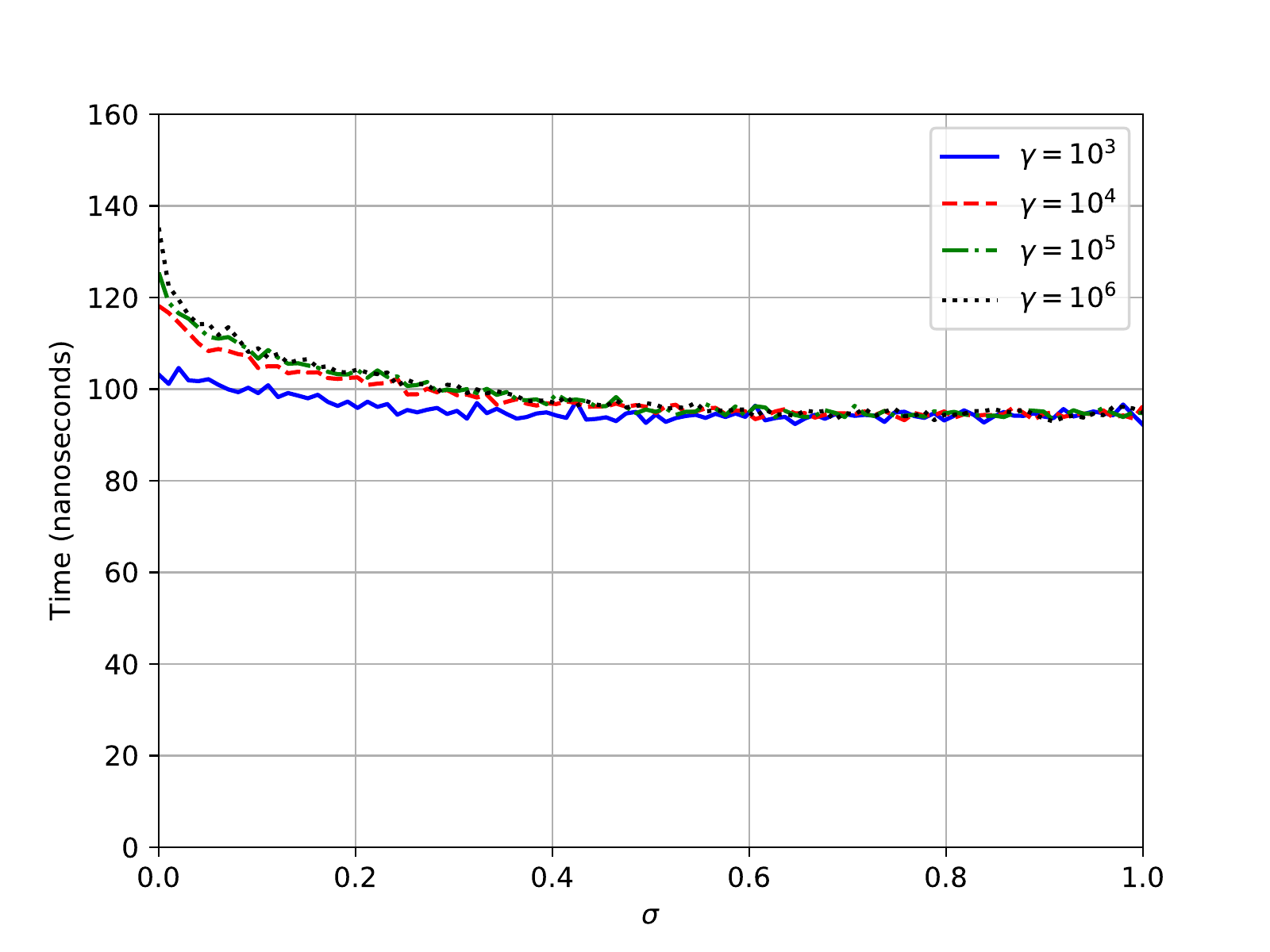}
\hfil

\vskip 3em

\hfil
\includegraphics[width=.49\textwidth]{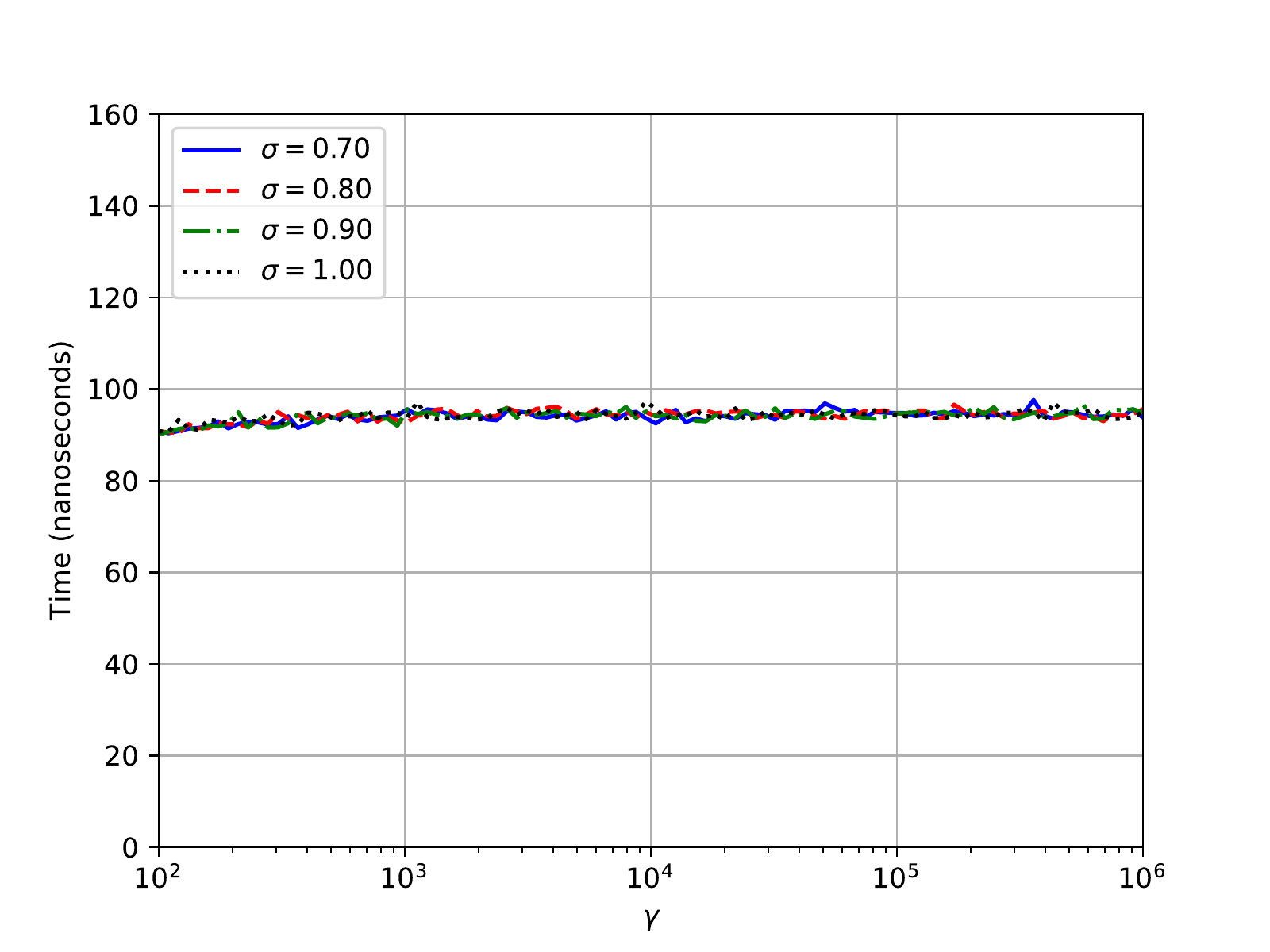}
\hfil
\includegraphics[width=.49\textwidth]{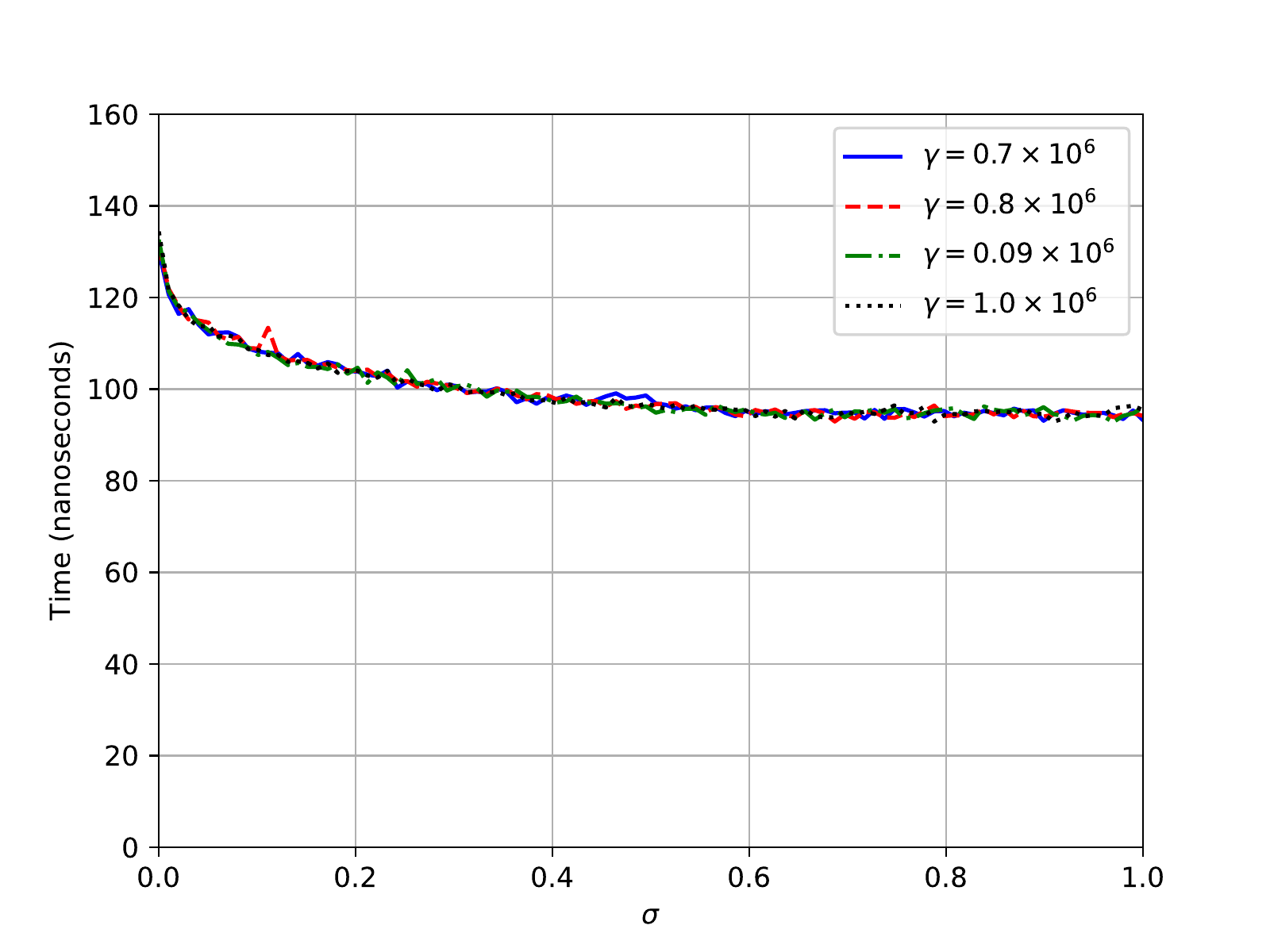}
\hfil

\caption{The time (in nanoseconds) required to evaluate  $\PS{n}(z;\gamma)$ using the algorithm
of this paper.
 Each of the graphs on the left gives the time required
as a function of $\gamma$ for several values of $\sigma$, while the plots on the right 
give the required time  as a function of $\sigma$ for several values of $\gamma$. 
A logarithmic scale is used for the x-axis in each of the plots on the left.}
\label{experiments:figure5}
\end{figure}


\end{document}